\documentclass{article}

\usepackage{latexsym}
\usepackage{amsmath}
\usepackage{amsfonts}
\usepackage{graphicx}

\newcommand{\ds}{\displaystyle}

\newcommand{\R}{\mathbb{R}}
\newcommand{\torus}{{\mathbb{T}}}
\newcommand{\gu}{\mathbf{u}}
\newcommand{\gm}{\mathbf{m}}

\usepackage{theorem}

\theoremheaderfont{\scshape}
\newtheorem{theorem}{Theorem}[section]

\theorembodyfont{\normalfont}
\newtheorem{remark}{Remark}[section]
\newenvironment{proof}{{\bf Proof }}{\hbox{~} \hfill \rule{0.5em}{0.5em}\\}
\numberwithin{equation}{section}

\begin{document}

\title{
Long term behaviour of singularly perturbed parabolic degenerated equation}         
\date{}          

\maketitle

\centerline{\scshape Ibrahima Faye \footnote{grandmbodj@hotmail.com }}
\medskip
{\footnotesize
 \centerline{Universit\'e de Bambey,UFR S.A.T.I.C, BP 30 Bambey (S\'en\'egal),}
\centerline{  Ecole Doctorale de 
                      Math\'ematiques et Informatique. }
   \centerline{Laboratoire de Math\'ematiques de la D\'ecision et d'Analyse Num\'erique}
\centerline{ (L.M.D.A.N) F.A.S.E.G)/F.S.T. }
} 

\medskip

\centerline{\scshape Emmanuel Fr\'enod \footnote{emmanuel.frenod@univ-ubs.fr}}
\medskip
{\footnotesize
 \centerline{Universit\'e Europ\'eenne de Bretagne, Lab-STICC (UMR CNRS 3192),}
 \centerline{ Universit\'e de Bretagne-Sud, Centre Yves Coppens,}
  \centerline{ Campus de Tohannic, F-56017, 
       Vannes Cedex, France}
  \centerline{ET}
   \centerline{Projet INRIA Calvi, Universit\'{e} de Strasbourg, IRMA,}
   \centerline{7 rue Ren\'e Descartes, F-67084 Strasbourg Cedex, France}
}

\medskip
\centerline{\scshape Diaraf SECK \footnote{diaraf.seck@ucad.edu.sn}}
\medskip
{\footnotesize
 \centerline{Universit\'e Cheikn Anta Diop de Dakar, BP 16889 Dakar Fann,}
  \centerline{  Ecole Doctorale de 
                      Math\'ematiques et Informatique. }
   \centerline{Laboratoire de Math\'ematiques de la D\'ecision et d'Analyse Num\'erique}
    \centerline{ (L.M.D.A.N) F.A.S.E.G)/F.S.T. }
      \centerline{ET}
   \centerline{UMMISCO, UMI 209, IRD, France}

}
\pagestyle{myheadings}
 \renewcommand{\sectionmark}[1]{\markboth{#1}{}}
\renewcommand{\sectionmark}[1]{\markright{\thesection\ #1}}
\begin{abstract}\noindent
In this paper we consider models built in \cite{FaFreSe} for short-term, mean-term 
and long-term morphodynamics of dunes and megariples. We give an existence and uniqueness result for long term dynamics of dunes.
This result is based on a time-space periodic solution existence result for degenerated parabolic equation that we set out. Finally the mean-term and long-term models are homogenized. \\
\end{abstract}
\section{Introduction and results}
In Faye, Fr\'enod and Seck \cite{FaFreSe}, based on works of Bagnold \cite{Bagnold}, Gadd, Lavelle and Swift \cite{GaddLavSw}, Idier\cite{Idier}, Astruc and Hulcher \cite{IdierAsH}, Meyer-Peter and Muller \cite{MePetMull} and Van Rijn \cite{Rijn1989}, we set out that a relevant model for short term dynamics of dunes, i.e. for their dynamics over several months, is
\begin{equation}\label{eq1}
   \frac{\partial z^\epsilon}{\partial t}-\frac{a}{\epsilon}\nabla\cdot\big((1-b\epsilon \gm)g_{a}(|\gu|)\nabla z^\epsilon\big)=\frac{c}{\epsilon}\nabla\cdot\left((1-b\epsilon \gm)g_{c}(|\gu|)\frac{\gu}{|\gu|}\right),
\end{equation}
where $z^{\epsilon} = z^{\epsilon}(x, t),$ is the dimensionless seabed altitude at $t$ and in  $x.$ For a given constant $T,\,\,t\in[0,T),$ 
stands for the dimensionless time and $x=(x_{1},x_{2})\in\torus^{2},$ $\torus^{2}$ being the two dimensional torus $\,\mathbb{R}^{2}/\,\mathbb{Z}^{2},$ is the dimensionless position variable. Operators $\nabla$ and $\nabla\cdot$ refer to gradient and divergence. Functions $g_{a}$ and $g_{c}$ are regular on $\R^+$ and satisfy
\begin{equation}\label{eq2}\left\{ \begin{array}{ccc}
g_{a}\geq g_{c}\geq0,\,\,g_{c}(0)=g'_{c}(0)=0,\\
\exists d\geq0,\sup_{u\in\mathbb{R}^{+}}|g_{a}(u)|+\sup_{u\in\mathbb{R}^{+}}|g'_{a}(u)|\leq d,\\
    \sup_{u\in\mathbb{R}^{+}}|g_{c}(u)|+\sup_{u\in\mathbb{R}^{+}}|g'_{c}(u)|\leq d,\\
   \exists U_{thr}\geq0,\,\,\exists G_{thr}>0,\,\,\textrm{such that}\,\, u\geq U_{thr}\Longrightarrow g_{a}(u)\geq G_{thr}.\end{array}\right.
\end{equation}

Fields $\gu$ and $\gm$ are dimensionless water velocity and height. They are given by
\begin{equation}\label{eq3}
\gu(t,x)=\mathcal{U}(t,\frac{t}{\epsilon},x),\quad \gm(t,x)=\mathcal{M}(t,\frac{t}{\epsilon},x),
 \end{equation}
where
\begin{equation}\label{eq4}\left\{ \begin{array}{ccc}
        \ds \mathcal{U}=\mathcal{U}(t,\theta,x)\,\,\textrm {and} \,\,\mathcal{M}=\mathcal{M}(t,\theta,x) \,\,\textrm{are regular functions on}\,\, \mathbb{R}^+\times\mathbb{R}\times\torus^{2},\\
        \ds\theta\longmapsto(\mathcal{U},\mathcal{M})\,\,\textrm{is periodic of period}\,\, 1,\\
        \ds |\mathcal{U}|,\,\,|\frac{\partial \mathcal{U}}{\partial t}|,\,\,|\frac{\partial \mathcal{U}}{\partial \theta}|,\,\,|\nabla \mathcal{U}|,\,\,|\mathcal{M}|,\,\,|\frac{\partial \mathcal{M}}{\partial t}|,\,\,|\frac{\partial \mathcal{M}}{\partial \theta}|,\,\,|\nabla \mathcal{M}| \,\,\textrm{are bounded by}\,\,d,\\
           \ds \forall (t,\theta,x)\in\mathbb{R}^{+}\times\mathbb{R}\times\torus^{2},\,\,|\mathcal{U}(t,\theta,x)|\leq U_{thr} \Longrightarrow
\hspace{3cm }\\ \hspace{3cm} \ds \frac{\partial\mathcal{U}}{\partial t}=0,\,\,\frac{\partial\mathcal{M}}{\partial t}=0,\,\,\nabla\mathcal{M}(t,\theta,x)=0\,\,\textrm{and}\,\,\nabla\mathcal{U}(t,\theta,x)=0,\\
            \ds\exists \theta_{\alpha}<\theta_{\omega}\in[0,1]\,\,\textrm{such that}\,\, \forall\,\,\theta\in [\theta_{\alpha},\theta_{\omega}]\Longrightarrow|\mathcal{U}(t,\theta,x)|\geq U_{thr} .\end{array}\right.
\end{equation}

\medskip 
A relevant model for mean term, i.e. when dune dynamics is observed over a few years, is
 
\begin{equation}\label{eq5}
\frac{\partial z^\epsilon}{\partial t}-\frac{a}{\epsilon}\nabla\cdot\left((1-b\sqrt{\epsilon}\gm)g_{a}(|\gu|)\nabla z^\epsilon\right)=\frac{c}{\epsilon}\nabla\cdot\left((1-b\sqrt{\epsilon}\gm)g_{c}(|\gu|)\frac{\gu}{|\gu|}\right),
\end{equation}
with condition (\ref{eq2}) on $g_{a}$ and $g_{c}$ and with  $\gu$ and $\gm$ given by
 \begin{equation}\label{eq6}
 \gu(t,x)=\widetilde{\mathcal{U}}(t,\frac{t}{\sqrt{\epsilon}},\frac{t}{\epsilon},x),\quad \gm(t,x)=\mathcal{M}(t,\frac{t}{\sqrt{\epsilon}},\frac{t}{\epsilon},x),
 \end{equation}
For mathematical reasons, we assumed
 \begin{equation}\label{eq7}
 \widetilde{\mathcal{U}}(t,\tau,\theta,x)=\mathcal{U}(t,\theta,x)+ \sqrt{\epsilon}\,\mathcal{U}_{1}(t,\tau,\theta,x),
 \end{equation}
where $\mathcal{U}=\mathcal{U}(t,\theta,x)$ and $\mathcal{U}_{1}=\mathcal{U}_{1}(t,\tau,\theta,x)$ are regular. We also assumed  that 
$\mathcal{M}=\mathcal{M}(t,\tau,\theta,x)$ is regular and
\begin{equation}\label{eq8}\left\{\begin{array}{ccc}
\ds\theta\longmapsto(\mathcal{U},\mathcal{U}_{1},\mathcal{M})\,\,\textrm{is periodic of period 1,}\\
\ds\tau\longmapsto(\mathcal{U}_{1},\mathcal{M})\,\,\textrm{is periodic of period 1,}\\
\ds |\mathcal{U}|,\,\,|\frac{\partial \mathcal{U}}{\partial t}|,\,\,|\frac{\partial \mathcal{U}}{\partial \theta}|,\,\,|\nabla \mathcal{U}|,\,\,|\mathcal{U}_{1}|,\,\,|\frac{\partial \mathcal{U}_{1}}{\partial t}|,\,\,|\frac{\partial \mathcal{U}_{1}}{\partial \tau}|,\,\,|\frac{\partial \mathcal{U}_{1}}{\partial \theta}|,\,\,|\nabla \mathcal{U}_{1}|,
\hspace{3cm}\\ \hspace{3cm}
\ds |\mathcal{M}|,\,\,|\frac{\partial \mathcal{M}}{\partial t}|,\,\,|\frac{\partial \mathcal{M}}{\partial \theta}|,\,\,|\frac{\partial \mathcal{M}}{\partial \tau}|,\,\,|\nabla \mathcal{M}| \,\,\textrm{are bounded by}\,\,d,\\
            \ds\forall (t,\tau,\theta,x)\in\mathbb{R}^{+}\times\mathbb{R}\times\mathbb{R}\times\torus^{2},\,\,|\widetilde{\mathcal{U}}(t,\tau,\theta,x)|\leq U_{thr}\Longrightarrow
 \hspace{2.5cm}\\  \hspace{0.5cm}
 \ds\frac{\partial\widetilde{\mathcal{U}}}{\partial t} (t,\tau,\theta,x) =0,\,\,\frac{\partial\widetilde{\mathcal{U}}}{\partial \tau} (t,\tau,\theta,x) =0,\,\,\nabla\widetilde{\mathcal{U}}(t,\tau,\theta,x)=0,
 \\ \hspace{2.5cm}
 \ds\frac{\partial\mathcal{M}}{\partial t} (t,\tau,\theta,x) =0,\,\,\frac{\partial\mathcal{M}}{\partial \tau} (t,\tau,\theta,x) =0\textrm{ and}\,\,\nabla\mathcal{M}(t,\tau,\theta,x)=0,\\
            \exists \theta_{\alpha}<\theta_{\omega}\in[0,1]\,\,\textrm{such that}\,\, \forall\,\,\theta\in [\theta_{\alpha},\theta_{\omega}]\Longrightarrow|\widetilde{\mathcal{U}}(t,\tau,\theta,x)|\geq U_{thr}.
\end{array}\right.\end{equation}

~

A relevant model for long-term
 dune dynamics is the following equation
\begin{equation}\label{eq81}
\frac{\partial z^\epsilon}{\partial t}-\frac{a}{\epsilon^{2}}\nabla\cdot\left((1-b\epsilon \gm)g_{a}(|\gu|)\nabla z^\epsilon\right)=\frac{c}{\epsilon^{2}}\nabla\cdot\left((1-b\epsilon \gm)g_{c}(|\gu|)\frac{\gu}{|\gu|}\right),
\end{equation}
where $a,$ $b$ and $c$ are constants, where $g_{a}$ and $g_{c}$ satisfy assumption (\ref{eq2}), and where $z^{\epsilon}$ is defined on the same space as before. It is also relevant to assume
\begin{equation}\label{eq9}
\gu(x,t)=\mathcal{U}(t,\frac{t}{\epsilon},x)=\mathcal{U}_{0}(\frac{t}{\epsilon})+\epsilon\mathcal{U}_{1}(\frac{t}{\epsilon},x)+\epsilon^{2}\mathcal{U}_{2}(t,\frac{t}{\epsilon},x),$$
$$\quad \gm(t,x)=\mathcal{M}(\frac{t}{\epsilon},x)+\epsilon^{2}\mathcal{M}_{2}(t,\frac{t}{\epsilon},x),
\end{equation}
where $\mathcal{U}_{0}=\mathcal{U}_{0}(\theta),\,\,\mathcal{U}_{1}=\mathcal{U}_{1}(\theta,x),\,\,\mathcal{U}_{2}=\mathcal{U}_{2}(t,\theta,x),\,\,\mathcal{M}=\mathcal{M}(\theta,x)$ and $\mathcal{M}_{2}=\mathcal{M}_{2}(t,\theta,x)$ are regular and
\begin{equation}\label{eq10}\left\{\begin{array}{ccc}
\theta\longmapsto(\mathcal{U}_0,\mathcal{U}_1,\mathcal{U}_{2},\mathcal{M},\mathcal{M}_{2})\,\,\textrm{is periodic of period 1},\\
\ds|\mathcal{U}_0|,\,\,|\frac{\partial \mathcal{U}_0}{\partial \theta}|,\,\,\ds|\mathcal{U}_1|,\,\,|\frac{\partial \mathcal{U}_1}{\partial \theta}|,\,\,|\nabla \mathcal{U}_1|,\,\,|\mathcal{U}_{2}|,\,\,|\frac{\partial \mathcal{U}_{2}}{\partial t}|,\,\,|\frac{\partial \mathcal{U}_{2}}{\partial \theta}|,\,\,|\nabla \mathcal{U}_{2}|,\,\,
\ds |\mathcal{M}|,\,\,|\frac{\partial \mathcal{M}}{\partial \theta}|,\,\,
\hspace{0.5cm}\\ 
\ds |\nabla \mathcal{M}|, \,\,  |\mathcal{M}_{2}|,\,\,|\frac{\partial \mathcal{M}_{2}}{\partial t}|,\,\,|\frac{\partial \mathcal{M}_{2}}{\partial \theta}|,\,\,|\nabla \mathcal{M}_{2}|\textrm{ are bounded by}\,\,d,\\
\ds \forall (t,\theta,x)\in\mathbb{R}^{+}\times\mathbb{R}\times\torus^{2},\,\,|\mathcal{U}_0(\theta)+\epsilon\mathcal{U}_1(\theta,x)+\epsilon^{2}\mathcal{U}_{2}(t,\theta,x)|\leq U_{thr}\Longrightarrow
\hspace{2cm}\\ \hspace{0.5cm}
\ds \frac{\partial{\mathcal{U}_{2}}}{\partial t} (t,\theta,x) =0,\,\,\nabla{\mathcal{U}_1}(\theta,x)=0,\,\,\nabla{\mathcal{U}_{2}}(t,\theta,x)=0,
\\ \hspace{3cm}
\ds\frac{\partial\mathcal{M}_{2}}{\partial t}(t,\theta,x)=0,\,\,\nabla\mathcal{M}(\theta,x)=0,\,\,\nabla\mathcal{M}_{2}(t,\theta,x)=0,
\\
\ds\exists \theta_{\alpha}<\theta_{\omega}\in[0,1]\,\,\textrm{such that}\,\, \forall\,\,\theta\in\mathbb{R},\,\,\theta\in [\theta_{\alpha},\theta_{\omega}] 
\hspace{5cm}\\ \hspace{5cm}
 \Longrightarrow|\mathcal{U}_0(\theta)+\mathcal{U}_1(\theta,x)+\epsilon^{2}\mathcal{U}_{2}(t,\theta,x)|\geq U_{thr}.
\end{array}\right.\end{equation}

~

Equations (\ref{eq1}), (\ref{eq5}) or (\ref{eq81}) need to be provided with an initial condition
\begin{equation}\label{eq11}
z^{\epsilon}_{|t=0}=z_{0},
\end{equation}
giving the shape of the seabed at the initial time.\\
In \cite{FaFreSe}, we then gave an existence and uniqueness result for short-term model (\ref{eq1}) if hypotheses (\ref{eq2}), (\ref{eq3}) and (\ref{eq4}) are satisfied and for the mean term one (\ref{eq5}), if hypotheses (\ref{eq2}), (\ref{eq6}), (\ref{eq7}) and (\ref{eq8}) are satisfied. This result was built on a time-space periodic solution existence result for degenerated parabolic equation. Under the same assumptions, the asymptotic behaviour of $z^{\epsilon},$ as $\epsilon\rightarrow0,$ solution of short term model (\ref{eq1}) is also given by homogenization methods.   
Futhermore if moreover  $U_{thr} = 0,$ a corrector result was set out, which gives a rigorous version of asymptotic expansion of the sequence $z^\epsilon$:
\begin{equation}
z^{\epsilon}(t,x)\,\,=\,\,U(t,\frac{t}{\epsilon},x)\,\,+\,\,\epsilon U^{1}(t,\frac{t}{\epsilon},x)+\dots, \;
\end{equation}
where $U$ and $U^{1}$ are solutions to
\begin{equation}\label{eq5.15}
\frac{\partial U}{\partial \theta}-\nabla\cdot\left(\widetilde{\mathcal{A}}\nabla U\right)=\nabla\cdot\widetilde{\mathcal{C}},
\end{equation}
\begin{equation}\label{eq5.5}
\frac{\partial U^{1}}{\partial \theta}-\nabla\cdot\left(\widetilde{\mathcal{A}}\nabla U^{1}\right)=\nabla\cdot\widetilde{\mathcal{C}}_{1}+\frac{\partial U}{\partial t}+\nabla\cdot(\widetilde{\mathcal{A}}_{1}\nabla U),
\end{equation}
where $\widetilde{\mathcal{A}}$ and $\widetilde{\mathcal{C}}$ are given by 
\begin{equation}
\widetilde{\mathcal{A}}=a\,g_a(|\mathcal{U}(t,\theta,x)|) \,\, \textrm{and}\,\,\widetilde{\mathcal{C}}=c\,g_c(|\mathcal{U}(t,\theta,x)|)\,\frac{\mathcal{U}(t,\theta,x)}{|\mathcal{U}(t,\theta,x)|},
\end{equation}

and  $\widetilde{\mathcal{A}}_{1}$ and $\widetilde{\mathcal{C}}_{1}$ are given by
\begin{multline}\widetilde{\mathcal{A}}_{1}(t,\theta,x)=-ab\mathcal{M}(t,\theta,x)\,g_a(|\mathcal{U}(t,\theta,x)|),\,\, 
\\ \text{ and }
\widetilde{\mathcal{C}}_{1}(t,\theta,x)=-cb\mathcal{M}(t,\theta,x)\,g_c(|\mathcal{U}(t,\theta,x)|)\,\frac{\mathcal{U}(t,\theta,x)}{|U(t,\theta,x)|}.\end{multline}
In \cite{FaFreSe}, we did not state neither any existence result for long term model (\ref{eq81}) nor any asymptotic behaviour result for mean term and long term models. Stating those results is the subject of the present paper. We will now state those main results. The first one concerns existence and
uniqueness for the long-term model. 
\begin{theorem}\label{th01}
For any $T>0,$ any $a>0,$ any real constants $b$ and $c$ and any $\epsilon>0,$ under assumptions (\ref{eq2}), (\ref{eq9}) and (\ref{eq10}), 
if
\begin{equation}\label{eq12}
z_{0}\,\,\in L^{2}(\torus^{2}),
\end{equation}
there exists a unique function $z^{\epsilon}\in L^{\infty}([0,T),L^{2}(\torus^{2})),$ solution to equation (\ref{eq81}) 
provided with initial condition (\ref{eq11}).\\
Moreover, for any $t\in[0,T],\,\, z^{\epsilon}$ satisfies
\begin{equation}\label{eq12.111}
\|z^{\epsilon}\|_{  L^{\infty}([0,T),L^{2}(\torus^{2}))}\leq\widetilde{\gamma},
\end{equation}
for a constant $\widetilde{\gamma}$ not depending on $\epsilon$ and
\begin{equation}\label{eq12.1}
\frac{\ds d\left(\int_{\torus^{2}} z^{\epsilon}(t,x) \, dx\right)}{dt} =0.
\end{equation}
\end{theorem}

The proof of this theorem is done in section \ref{secExEs1}, except equality (\ref{eq12.1})  which is directly gotten by integrating (\ref{eq81})  with respect to $x$ over $\torus^2$.\\

\medskip 

We now give a result concerning the asymptotic behaviour as $\epsilon\rightarrow0$ of the
long term model. We notice that, since $\mathcal{U}$ and $\mathcal{M}+\epsilon^2\mathcal{M}_2$ do not depend on $t$ and $x$ when $\mathcal{U}\leq U_{thr},$ we have the following property:
\begin{equation}\label{suppl1} \forall \theta\in [0,1],\,\,\Big(\exists (t,x)\in[0,T)\times \torus^2\,\,\text{such that}\,\,\mathcal{U}(t,\theta,x)=0\,\,\text{or}\,\,\mathcal{M}(\theta,x)+\epsilon^2\mathcal{M}_2(t,\theta,x)=0\Big)$$
$$\Longrightarrow\Big(\forall (t,x)\in[0,T)\times \torus^2,\,\,\mathcal{U}(t,\theta,x)=0\,\,\text{and}\,\,\mathcal{M}(\theta,x)+\epsilon^2\mathcal{M}_2(t,\theta,x)=0\Big),
\end{equation}
and
\begin{equation}\label{suppl2}\{\theta\in [0,1],\,\,\mathcal{U}(\cdot,\theta,\cdot)=0\,\,\text{and}\,\,\mathcal{M}(\theta,\cdot)+\epsilon^2\mathcal{M}_2(\cdot,\theta,\cdot)=0\}$$
$$ \text{is an union of several intervals}.\end{equation}
Moreover we denote
\begin{equation}\label{suppl3} \Theta=[0,T)\times\{\theta\in\mathbb{R},\,\,\mathcal{U}(\cdot,\theta,\cdot)=0\,\,\text{and}\,\,\mathcal{M}(\theta,\cdot)+\epsilon^2\mathcal{M}_2(\cdot,\theta,\cdot)=0\}\times\torus^2,
\end{equation}
and
\begin{equation}\label{suppl3.1} \Theta_{thr}=\{(t,\theta,x)\in [0,T)\times\mathbb{R}\times\torus^2,\,\,\mathcal{U}(t,\theta,x)<U_{thr}\}.
\end{equation}

\begin{theorem}\label{thAsyBeh1}
For any $T>0,$ under the same assumptions as in theorem \ref{th01}, the sequence of solutions $(z^\epsilon)$ to
equation (\ref{eq81}) given by theorem \ref{th01} two-scale converges to a profile\\
$U\in L^{\infty}([0,T],L^{\infty}_\#(\R,L^2(\torus^{2})))$
 which is the unique solution to
\begin{equation}\label{eqHomIntro1}
-\nabla\cdot(\widetilde{\mathcal{A}}\nabla U)=\nabla \cdot\widetilde{\mathcal{C}}\,\,\text{on}\,\,\Big([0,T)\times\mathbb{R}\times\torus^2\Big)\backslash\Theta,
\end{equation}
\begin{equation}\label{eqHomIntroa}
\frac{\partial U}{\partial\theta}\,\,=\,\,0\,\,\text{on}\,\, \Theta_{thr},
\end{equation}
and\begin{equation}\int_0^1\int_{\torus^2}U\,d\theta\,dx=\int_{\torus^2}z_0dx,
\end{equation}
where $\widetilde{\mathcal{A}}$ and $\widetilde{\mathcal{C}}$ are given by
\begin{equation}\label{eqHomIntroBis}
\widetilde{\mathcal{A}}=a\,g_a(|\mathcal{U}(t,\theta,x)|) \,\, \textrm{and}\,\,\widetilde{\mathcal{C}}=c\,g_c(|\mathcal{U}(t,\theta,x)|)\,\frac{\mathcal{U}(t,\theta,x)}{|\mathcal{U}(t,\theta,x)|}.
\end{equation}
\end{theorem}

Above and in the sequel, for all $p\geq1$ and $q\geq1,$ we denote by\\ $L^p_\#(\mathbb{R}, L^q(\torus^2))=\Big\{f:\mathbb{R}\longrightarrow L^q(\torus^2)\,\,\text{mesurable and periodic of period 1 in $\theta$}\,\,\text{such that}\,\,\\\theta\mapsto\|f(\theta)\|_{L^{q}(\torus^2)}\in L^{p}([0,1])\Big\}.$
\begin{remark}
Notice that $\Big([0, T)\times\mathbb{R}\times\torus^2\Big)\backslash\Theta\cap\Theta_{thr}$ is not empty. On this set $0<\mathcal{U}< U_{thr}.$

\medskip 
 This contributes to make of (\ref{eqHomIntro1}),(\ref{eqHomIntroa}) a well posed problem.\\ 

\end{remark}

\medskip

\noindent Now we turn to mean term model for which we set out asymptotic behaviours.
\begin{theorem}\label{thHomSecHom1}
Under assumptions (\ref{eq2}), (\ref{eq6}), (\ref{eq7}) and (\ref{eq8}), for any $T,$ not depending on $\epsilon,$  the sequence $(z^{\epsilon})$ of solutions to (\ref{eq5}) built in \cite{FaFreSe} provided with initial condition (\ref{eq11})
two-scale converges to the profile $U\in L^{\infty}([0,T]\times\mathbb{R},L^{\infty}_\#(\R,L^2(\torus^{2})))$ solution to
\begin{equation}
\label{ee179b}
\frac{\partial U}{\partial\theta}-\nabla\cdot(\widetilde{\mathcal{A}}\nabla U)=\nabla \cdot\widetilde{\mathcal{C}},
\end{equation} where $\widetilde{\mathcal{A}}$ and $\widetilde{\mathcal{C}}$ are given by 
\begin{equation}\label{eqHomIntroBis2}
\widetilde{\mathcal{A}}=a\,g_a(|\mathcal{U}(t,\tau,\theta,x)|) \,\, \textrm{and}\,\,\widetilde{\mathcal{C}}=c\,g_c(|\mathcal{U}(t,\tau,\theta,x)|)\,\frac{\mathcal{U}(t,\tau,\theta,x)}{|\mathcal{U}(t,\tau,\theta,x)|}.
\end{equation}
 
\end{theorem}
Finally, a corrector result for the mean-term model is given under restrictive assumptions.

\begin{theorem}\label{th3.1}
Under assumptions (\ref{eq2}), (\ref{eq6}), (\ref{eq7}), (\ref{eq8}) and if moreover $U_{thr}=0,$ considering 
function $z^{\epsilon}\in L^{\infty}([0,T),L^{2}(\torus^{2})),$ 
solution to (\ref{eq5}) with initial condition (\ref{eq11}) and function $U^{\epsilon}\in L^{\infty}([0,T),L^{2}(\torus^{2}))$ defined by 
\begin{equation}
U^{\epsilon}(t,x)=U(t,\frac{t}{\sqrt{\epsilon}},\frac{t}{\epsilon},x),
\end{equation} 
where $U$ is the solution to (\ref{ee179b}), the following estimate is satisfied:
\begin{equation}
\Big\|\frac{z^{\epsilon}-U^{\epsilon}}{\sqrt{\epsilon}}\Big\|_{  L^{\infty}([0,T),L^{2}(\torus^{2}))}\,\,\leq\alpha,
\end{equation}
where $\alpha$ is a constant not depending on $\epsilon.$\\
Furthermore, 
\begin{equation}
\frac{z^{\epsilon}-U^{\epsilon}}{\sqrt{\epsilon}}\quad\textrm{ two-scale converges to a profile}\,\,U_{\frac{1}{2}}\in L^{\infty}([0,T]\times\mathbb{R},L^{\infty}_\#(\R,L^2(\torus^{2}))),
\end{equation} 
which is the unique solution to
\begin{equation}\frac{\partial
U_{\frac{1}{2}}}{\partial\theta}-\nabla\cdot\Big(\widetilde{\mathcal{A}}\nabla
U_{\frac{1}{2}}\Big)=\nabla\cdot\widetilde{\mathcal{C}}_{1}\\
+\nabla\cdot\Big(\widetilde{\mathcal{A}}_{1}\nabla
U\Big)-\frac{\partial U}{\partial\tau}\end{equation}
where $\widetilde{\mathcal{A}}$ and $\widetilde{\mathcal{C}}$ are given by (\ref{eqHomIntroBis2}) and where $\widetilde{\mathcal{A}}_{1}$ and $\widetilde{\mathcal{C}}_{1}$ are given by
\begin{multline}\widetilde{\mathcal{A}}_{1}(t,\tau,\theta,x)=-ab\mathcal{M}(t,\tau,\theta,x)\,g_a(|\mathcal{U}(t,\theta,\tau,x)|),\,\, 
\\ \text{ and }
\widetilde{\mathcal{C}}_{1}(t,\tau,\theta,x)=-cb\mathcal{M}(t,\tau,\theta,x)\,g_c(|\mathcal{U}(t,\tau,\theta,x)|)\,\frac{\mathcal{U}(t,\tau,\theta,x)}{|U(t,\tau,\theta,x)|}.\end{multline}
\end{theorem}

\section{Existence and estimates,\,\,proof of theorem \ref{th01}} \label{secExEs1}      
Setting:
\begin{equation}\label{La1}
    \mathcal{A}^{\epsilon}(t,x)=\widetilde{\mathcal{A}}_{\epsilon}(t,\frac{t}{\epsilon},x),
\end{equation}
and \begin{equation}\label{La4}
    \mathcal{C}^{\epsilon}(t,x)=\widetilde{\mathcal{C}}_{\epsilon}(t,\frac{t}{\epsilon},x),
\end{equation}
where
\begin{equation}\label{La3}\widetilde{\mathcal{A}}_{\epsilon}(t,\theta,x)
=a(1-b\epsilon\mathcal{M}(t,\theta,x))\,
g_a(|\mathcal{U}(t,\theta,x)|),
\end{equation} 
and
\begin{equation}\label{La5} \widetilde{\mathcal{C}}_{\epsilon}(t,\theta,x)=
c(1-b\epsilon\mathcal{M}(t,\theta,x))\,g_c(|\mathcal{U}(t,\theta,x)|)\,
\frac{\mathcal{U}(t,\theta,x)}{|\mathcal{U}(t,\theta,x)|},
\end{equation}
equation (\ref{eq81}) with initial condition (\ref{eq11}) can be set in the form
\begin{equation}\label{L1}\left\{\begin{array}{cc}
    \ds \frac{\partial z^{\epsilon}}{\partial t}-\frac{1}{\epsilon^{2}}\nabla\cdot(\mathcal{A}^{\epsilon}\nabla z^{\epsilon})=\frac{1}{\epsilon^{2}}\nabla\cdot\mathcal{C}^{\epsilon},\\
    \ds z^{\epsilon}_{|t=0}=z_{0}.
\end {array}\right.\end{equation}
Because of hyptothesis (\ref{eq81}) and under assumptions (\ref{eq2}) and (\ref{eq10}),  $\widetilde{\mathcal{A}}_{\epsilon}$ and $\widetilde{\mathcal{C}}_{\epsilon}$
given by (\ref{La3}) and (\ref{La5})  satisfy the following hypotheses:
\begin{equation}\label{Hy1}\left\{\begin{array}{ccc}
\theta\longmapsto(\widetilde{\mathcal{A}}_{\epsilon},\widetilde{\mathcal{C}}_{\epsilon})\,\,\textrm{ is periodic of period}\,\,1,
\\

x\longmapsto(\widetilde{\mathcal{A}}_{\epsilon},\widetilde{\mathcal{C}}_{\epsilon})\,\,\textrm{is defined on}\,\,
\torus^2,\,\,\\

\ds|\widetilde{\mathcal{A}}_{\epsilon}|\leq\gamma,\,\,|\widetilde{\mathcal{C}}_{\epsilon}|\leq\gamma,\,\,\left|\frac{\partial\widetilde{\mathcal{A}}_{\epsilon}}{\partial t}\right|\leq\epsilon^2\gamma,\,\,\left|\frac{\partial\widetilde{\mathcal{C}}_{\epsilon}}{\partial t}\right|\leq\epsilon^2\gamma,\left|\frac{\partial\nabla\widetilde{\mathcal{A}}_{\epsilon}}{\partial t}\right|\leq\epsilon^2\gamma,\\

\ds\left|\frac{\partial\widetilde{\mathcal{A}}_{\epsilon}}{\partial \theta}\right|\leq\gamma,\,\,
    \left|\frac{\partial\widetilde{\mathcal{C}}_{\epsilon}}{\partial \theta}\right|\leq\gamma,\,\,|\nabla\widetilde{\mathcal{A}}_{\epsilon}|\leq\epsilon\gamma,\,\,|\nabla\cdot\widetilde{\mathcal{C}}_{\epsilon}|\leq\epsilon\gamma,\,\,\left|\frac{\partial\nabla\cdot\widetilde{\mathcal{C}}_{\epsilon}}{\partial t}\right|\leq\epsilon^2\gamma,\\
\end{array}\right.\end{equation}
\begin{equation}\label{Hy2}\left\{\begin{array}{ccc}
\exists \widetilde{G}_{thr}, ~ \theta_{\alpha}<\theta_{\omega}\in [0,1]\,\,\textrm{such that}\,\, 
\forall\,\,\theta\in [\theta_{\alpha},\theta_{\omega}]\Longrightarrow
 \widetilde{\mathcal{A}}_{\epsilon}(t,\theta,x)\geq\widetilde{G}_{thr},
 \\
 \ds
 \widetilde{\mathcal{A}}_{\epsilon}(t,\theta,x)\leq\widetilde{G}_{thr}\Longrightarrow\left\{\begin{array}{ccc}
\ds\frac{\partial\widetilde{\mathcal{A}}_{\epsilon}}{\partial t}(t,\theta,x)=0,\,\,
\nabla\widetilde{\mathcal{A}}_{\epsilon}(t,\theta,x)=0, \vspace{3pt}\\
\ds\frac{\partial\widetilde{\mathcal{C}}_{\epsilon}}{\partial t}(t,\theta,x)=0,\,\,
\nabla\cdot\widetilde{\mathcal{C}}_{\epsilon}(t,\theta,x)=0 ,\\
\end{array}\right.\end{array}\right.\end{equation}
and
 \begin{equation}\label{Hy3}\left\{\begin{array}{ccc}
\ds
 |\widetilde{\mathcal{C}}_{\epsilon}|\leq\gamma|\widetilde{\mathcal{A}}_{\epsilon}|, ~
 |\widetilde{\mathcal{C}}_{\epsilon}|^{2}\leq\gamma|\widetilde{\mathcal{A}}_{\epsilon}|,~
 |\nabla\widetilde{\mathcal{A}}_{\epsilon}|\leq\epsilon\gamma|\widetilde{\mathcal{A}}_{\epsilon}|,~
 \Big|\frac{\partial\widetilde{\mathcal{A}}_{\epsilon}}{\partial t}\Big|\leq\epsilon^2\gamma|\widetilde{\mathcal{A}}_{\epsilon}|,
\\
\ds
 \Big|\frac{\partial(\nabla\widetilde{\mathcal{A}}_{\epsilon})}{\partial t}
       \Big|^{2}\leq\epsilon^2\gamma|\widetilde{\mathcal{A}}_{\epsilon}|,~
 \Big|\nabla \cdot\widetilde{\mathcal{C}}_{\epsilon}\Big|\leq\epsilon\gamma|\widetilde{\mathcal{A}}_{\epsilon}|,~
 \Big|\frac{\partial\widetilde{\mathcal{C}}_{\epsilon}}{\partial t}\Big|\leq\epsilon^2\gamma|\widetilde{\mathcal{A}}_{\epsilon}|,~
 \Big|\frac{\partial\widetilde{\mathcal{C}}_{\epsilon}}{\partial t}\Big|^{2}\leq\epsilon^2\gamma^2|\widetilde{\mathcal{A}}_{\epsilon}|.
\end{array}\right.\end{equation}
 
\medskip

In this section we focus on existence and uniqueness of time-space periodic parabolic equations. From this, we then get existence of solution to equation (\ref{L1}). Existence of $z^{\epsilon}$ over a time interval depending on $\epsilon,$ is a traightforward consequence  of adaptations of results from LadyzensKaja, Solonnikov and Ural' Ceva \cite{LSU} or Lions \cite{J.L.L}.
Our aim is to proove that $z^{\epsilon}$ solution to (\ref{L1}) is bounded indepently of $\epsilon.$  We are going to introduce the following regularized equations. We recall that the method used is similar to the one used in \cite{FaFreSe}. 
\begin{equation}\label{L6} 
\frac{\partial \mathcal{S}^{\nu}}{\partial \theta}-\frac{1}{\epsilon}\nabla\cdot\Big(\big(\tilde{\mathcal{A}}_{\epsilon}(t,\cdot,\cdot)+\nu\big)\nabla\mathcal{S}^{\nu}\Big)=\frac{1}{\epsilon}\nabla\cdot\tilde{\mathcal{C}}_{\epsilon}(t,\cdot,\cdot),
\end{equation} and
\begin{equation} \label{L7}
\mu\mathcal{S}^{\nu}_{\mu}+\frac{\partial\mathcal{S}^{\nu}_{\mu}}{\partial\theta}-\frac{1}{\epsilon}\nabla\cdot\Big(\big(\tilde{\mathcal{A}}_{\epsilon}(t,\cdot,\cdot)+\nu\big)\nabla\mathcal{S}^{\nu}_{\mu}\Big)=\frac{1}{\epsilon}\nabla\cdot\tilde{\mathcal{C}}_{\epsilon}(t,\cdot,\cdot),
\end{equation}
where $\mu$ and $\nu$ are positive parameters.\\
We first prove existence of solutions $\mathcal{S}^{\nu}_{\mu}$  of (\ref{L7}) and we give estimates of $\mathcal{S}^{\nu}_{\mu}.$ 

\begin{theorem}\label{theoexes12}
Under assumptions (\ref{Hy1}), (\ref{Hy2}) and (\ref{Hy3}), for any $\mu>0$ and any $\nu>0,$ there exists a unique $\mathcal{S}_{\mu}^{\nu}=\mathcal{S}_{\mu}^{\nu}(t,\theta,x)\in\mathcal{C}^{0}\cap L^{2}(\mathbb{R}\times\torus^{2})$, periodic of period 1 with respect to $\theta,$  solution to (\ref{L7})
and regular with respect to the parameter $t.$ Moreover, the following estimates are satisfied
\begin{equation}\label{3.280}
    \sup_{\theta\in\mathbb{R}}\left|\int_{\torus^{2}}\mathcal{S}_{\mu}^{\nu}(\theta,x)dx\right|=0,
\end{equation}
\begin{equation}\label{+01}
 \|\nabla\mathcal{S}_{\mu}^{\nu}\|_{L^{2}_{\#}(\mathbb{R},L^{2}(\torus^{2}))}\leq\frac{\gamma}{\nu},
\end{equation}
\begin{equation}\label{+02}
 \|\Delta\mathcal{S}_{\mu}^{\nu}\|_{L^{2}_{\#}(\mathbb{R},L^{2}(\torus^{2}))}\leq\sqrt{2}\frac{\epsilon\gamma}{\nu}\sqrt{\frac{\gamma^2}{\nu^2}+1},
\end{equation}

\begin{equation}\label{3.281}
    \left\|\frac{\partial\mathcal{S}_{\mu}^{\nu}}{\partial \theta}\right\|_{L^{2}_{\#}(\mathbb{R},L^{2}(\torus^{2}))}\leq\frac{\gamma}{\sqrt{\epsilon\nu}}\sqrt{(\frac{\gamma}{2\nu}+1)},
\end{equation}

\begin{equation}\label{3.282}
    \|\nabla\mathcal{S}_{\mu}^{\nu}\|_{L^{\infty}_{\#}(\mathbb{R},L^{2}(\torus^{2}))}\leq\sqrt{\frac{\gamma^2}{\nu^2}+\frac{2\epsilon\gamma^2}{\nu}(\frac{\gamma^2}{\nu^2}+1)},
\end{equation}
\begin{equation}\label{3.283}
    \|\mathcal{S}_{\mu}^{\nu}\|_{L^{\infty}_{\#}(\mathbb{R},L^{2}(\torus^{2}))}\leq\sqrt{\frac{\gamma^2}{\nu^2}+\frac{2\epsilon\gamma^2}{\nu}(\frac{\gamma^2}{\nu^2}+1)},
\end{equation}
\begin{equation}\label{3.29}
    \left\|\frac{\partial\mathcal{S}_{\mu}^{\nu}}{\partial t}\right\|_{L^{\infty}_{\#}(\mathbb{R},L^{2}(\torus^{2}))}\leq\epsilon^3\frac{\gamma}{\nu}(1+\frac{\gamma}{\nu}).\end{equation}

\end{theorem}

\begin{proof}.\,(of Theorem \ref{theoexes12}).\,\,
The proof of this theorem is very similar to the one of Theorem 3.3 of Faye, Fr\'enod and Seck \cite{FaFreSe}. The big difference is the presence of $\frac{1}{\epsilon}-$ factors in (\ref{L7}). Hence we only sketch the most similar arguments and focus on the management of those $\frac{1}{\epsilon}-$factors.\\ 
Integrating equation (\ref{L7}) over $\mathbb{T}^{2}$ gives 
\begin{equation}
\mu\int_{\torus^{2}}\mathcal{S}_{\mu}^{\nu}dx+\int_{\torus^{2}}\frac{\partial\mathcal{S}^{\nu}_{\mu}}{\partial\theta}dx-\frac{1}{\epsilon}\int_{\torus^{2}}\nabla\cdot\Big((\tilde{\mathcal{A}}_{\epsilon}+\nu)\nabla\mathcal{S}^{\nu}_{\mu}\Big)dx=\frac{1}{\epsilon}\int_{\torus^{2}}\nabla\cdot\tilde{\mathcal{C}}_{\epsilon}dx,
\end{equation}
then $$\mu\int_{\torus^{2}}\mathcal{S}_{\mu}^{\nu}dx+\frac{d(\int_{\torus^{2}}\mathcal{S}_{\mu}^{\nu}dx)}{\partial\theta}=0, $$ which gives
$$\int_{\torus^{2}}\mathcal{S}_{\mu}^{\nu}(\theta,x)dx=\int_{\torus^{2}}\mathcal{S}_{\mu}^{\nu}(\tilde{\theta},x)e^{-\mu(\theta-\tilde{\theta})}dx.$$ Since $\mathcal{S}_{\mu}^{\nu}$ is periodic of period $1$ with respect to $\theta,$ $\int_{\torus^{2}}\mathcal{S}_{\mu}^{\nu}(\theta,x)dx $ is also periodic of period 1. Then (\ref{3.280}) is true.\\
Multiplying equation (\ref{L7}) by $\mathcal{S}_{\mu}^{\nu},$ integrating over  $\torus^{2}$  
and from $0$ to $1$ with respect to  $\theta$ gives
$$\mu\|\mathcal{S}_{\mu}^{\nu}\|_{L^{2}_{\#}(\mathbb{R},L^{2}(\torus^{2}))}^{2}+\frac{1}{2}\int_{0}^{1}\frac{d\|\mathcal{S}_{\mu}^{\nu}\|_{2}^{2}}{d\theta}d\theta+\frac{1}{\epsilon}\int_{0}^{1}\int_{\torus^{2}}(\tilde{\mathcal{A}}_{\epsilon}+\nu)|\nabla\mathcal{S}^{\nu}_{\mu}|^{2}dxd\theta\leq \frac{\gamma}{\epsilon}\int_{0}^{1}\int_{\torus^{2}}|\nabla\mathcal{S}^{\nu}_{\mu}|dxd\theta.$$
Since  $\tilde{\mathcal{A}}_{\epsilon}+\nu\geq\nu$ and taking into account that the above first term is positive and the second one equals zero, we have 
$$\frac{\nu}{\epsilon}\int_{0}^{1}\int_{\torus^{2}}|\nabla\mathcal{S}^{\nu}_{\mu}|^{2}dxd\theta\leq\frac{\gamma}{\epsilon}\|\nabla\mathcal{S}^{\nu}_{\mu}\|_{L^{2}_{\#}(\mathbb{R},L^{2}(\torus^{2}))},$$
then
$$\|\nabla\mathcal{S}^{\nu}_{\mu}\|_{L^{2}_{\#}(\mathbb{R},L^{2}(\torus^{2}))}^{2}\leq\frac{\gamma}{\nu}\|\nabla\mathcal{S}^{\nu}_{\mu}\|_{L^{2}_{\#}(\mathbb{R},L^{2}(\torus^{2}))},$$
which gives (\ref{+01}).\\
Multiplying  (\ref{L7}) by $\frac{\partial\mathcal{S}^{\nu}_{\mu}}{\partial\theta},$ integrating over $\torus^{2}$ and
integrating from $0$ to $1$ with respect to $\theta$ gives
\begin{gather}\Big\|\frac{\partial\mathcal{S}_{\mu}^{\nu}}{\partial\theta}\Big\|^{2}_{L^{2}_{\#}(\mathbb{R},L^{2}(\torus^{2}))}=\frac{1}{2\epsilon}\int_{0}^{1}\int_{\torus^{2}}\frac{\partial \tilde{\mathcal{A}}_{\epsilon}}{\partial\theta}|\nabla\mathcal{S}_{\mu}^{\nu}|^{2}dxd\theta+\frac{1}{\epsilon}\int_{0}^{1}\int_{\torus^{2}}\frac{\partial \tilde{\mathcal{C}}_{\epsilon}}{\partial\theta}\nabla\mathcal{S}_{\mu}^{\nu}dxd\theta\\
\leq \frac{\gamma}{\epsilon}\Big(\frac{1}{2}\|\nabla\mathcal{S}_{\mu}^{\nu}\|^{2}_{L^{2}_{\#}(\mathbb{R},L^{2}(\torus^{2}))}+\|\nabla\mathcal{S}_{\mu}^{\nu}\|_{L^{2}_{\#}(\mathbb{R},L^{2}(\torus^{2}))}\Big),
\end{gather}
which gives (\ref{3.281}).\\
Multiplying  (\ref{L7}) by $-\Delta\mathcal{S}_{\mu}^{\nu},$ and integrating over  $\torus^{2}$ gives
\begin{gather*}\mu\int_{\torus^{2}}|\nabla\mathcal{S}_{\mu}^{\nu}|^{2}dx+\int_{\torus^{2}}\nabla\mathcal{S}_{\mu}^{\nu}\cdot\nabla
\big(\frac{\partial\mathcal{S}_{\mu}^{\nu}}{\partial\theta}\big)dx+\frac{1}{\epsilon}\int_{\torus^{2}}\nabla\tilde{\mathcal{A}}_{\epsilon}\cdot\nabla\mathcal{S}_{\mu}^{\nu}\Delta\mathcal{S}_{\mu}^{\nu}dx+\\\frac{1}{\epsilon}\int_{\torus^{2}}(\tilde{\mathcal{A}}_{\epsilon}+\nu)
|\Delta\mathcal{S}_{\mu}^{\nu}|^{2}dx=-\frac{1}{\epsilon}\int_{\torus^{2}}\nabla\cdot\tilde{\mathcal{C}}_{\epsilon}\Delta\mathcal{S}_{\mu}^{\nu}dx,\end{gather*}
or
\begin{gather*}\mu\|\nabla\mathcal{S}_{\mu}^{\nu}\|_{2}^{2}+\frac{1}{2}\frac{d\|\nabla\mathcal{S}_{\mu}^{\nu}\|_{2}^{2}}{d\theta}+\frac{1}{\epsilon}\int_{\torus^{2}}(\tilde{\mathcal{A}}_{\epsilon}+\nu)
|\Delta\mathcal{S}_{\mu}^{\nu}|^{2}dx=\\-\frac{1}{\epsilon}\int_{\torus^{2}}\nabla\tilde{\mathcal{A}}_{\epsilon}\cdot\nabla\mathcal{S}_{\mu}^{\nu}\Delta\mathcal{S}_{\mu}^{\nu}dx -\frac{1}{\epsilon}\int_{\torus^{2}}\nabla\cdot\tilde{\mathcal{C}}_{\epsilon}\Delta\mathcal{S}_{\mu}^{\nu}dx.\end{gather*}
Since for any real number $U$ and $V$ 
\begin{equation}\label{eqa}
|UV|\leq\frac{\tilde{\mathcal{A}}_{\epsilon}+\nu}{4\epsilon}U^{2}+\frac{\epsilon}{\tilde{\mathcal{A}}_{\epsilon}+\nu}V^{2},
\end{equation}
using this formula with  $U=\Delta\mathcal{S}_{\mu}^{\nu},\,\,V=\frac{\nabla\tilde{\mathcal{A}}_{\epsilon}\cdot\nabla\mathcal{S}_{\mu}^{\nu}}{\epsilon},$ we have
$$\frac{1}{\epsilon}\int_{\torus^{2}}\nabla\tilde{\mathcal{A}}_{\epsilon}\cdot\nabla\mathcal{S}_{\mu}^{\nu}\Delta\mathcal{S}_{\mu}^{\nu}dx\leq\int_{\torus^{2}}\frac{\tilde{\mathcal{A}}_{\epsilon}+\nu}{4\epsilon}|\Delta\mathcal{S}_{\mu}^{\nu}|^{2}dx+\int_{\torus^{2}}\frac{1}{\epsilon(\tilde{\mathcal{A}}_{\epsilon}+\nu)} |\nabla\tilde{\mathcal{A}}_{\epsilon}\cdot\nabla\mathcal{S}_{\mu}^{\nu}|^{2}dx.$$
Taking   $U=\Delta\mathcal{S}_{\mu}^{\nu},\,\,V=\frac{\nabla\cdot\tilde{\mathcal{C}}_{\epsilon}}{\epsilon}$ and using again (\ref{eqa}) we obtain
$$\frac{1}{\epsilon}\int_{\torus^{2}}\nabla\cdot\tilde{\mathcal{C}}_{\epsilon}\Delta\mathcal{S}_{\mu}^{\nu}\leq\int_{\torus^{2}}\frac{\tilde{\mathcal{A}}_{\epsilon}+\nu}{4\epsilon}|\Delta\mathcal{S}_{\mu}^{\nu}|^{2}dx+\int_{\torus^{2}}\frac{1}{\epsilon(\tilde{\mathcal{A}}_{\epsilon}+\nu)} |\nabla\cdot\tilde{\mathcal{C}}_{\epsilon}|^{2}dx.$$
These two results give
\begin{equation}\mu\|\nabla\mathcal{S}_{\mu}^{\nu}\|_{2}^{2}+\frac{1}{2}\frac{d\|\nabla\mathcal{S}_{\mu}^{\nu}\|_{2}^{2}}{d\theta}+\frac{1}{\epsilon}\int_{\torus^{2}}(\tilde{\mathcal{A}}_{\epsilon}+\nu)
|\Delta\mathcal{S}_{\mu}^{\nu}|^{2}dx\leq   $$
$$\int_{\torus^{2}}\frac{\tilde{\mathcal{A}}_{\epsilon}+\nu}{2\epsilon}|\Delta\mathcal{S}_{\mu}^{\nu}|^{2}dx+\frac{1}{\epsilon}\int_{\torus^{2}}\frac{1}{\epsilon(\tilde{\mathcal{A}}_{\epsilon}+\nu)} \Big(|\nabla\tilde{\mathcal{A}}_{\epsilon}\cdot\nabla\mathcal{S}_{\mu}^{\nu}|^{2} +   |\nabla\cdot\tilde{\mathcal{C}}_{\epsilon}|^{2}\Big)dx,\end{equation}
or, using (\ref{Hy1}),
\begin{equation}\label{L15}\mu\|\nabla\mathcal{S}_{\mu}^{\nu}\|_{2}^{2}+\frac{1}{2}\frac{d\|\nabla\mathcal{S}_{\mu}^{\nu}\|_{2}^{2}}{d\theta}+\int_{\torus^{2}}\frac{(\tilde{\mathcal{A}}_{\epsilon}+\nu)}{2\epsilon}
|\Delta\mathcal{S}_{\mu}^{\nu}|^{2}dx\leq\frac{\epsilon^2\gamma^{2}}{\nu\epsilon}\Big(\int_{\torus^{2}}|\nabla\mathcal{S}_{\mu}^{\nu}|^{2}dx+1\Big),\end{equation}
and integrating from  $0$ to $1$ with respect to $\theta,$ we have
$$\mu\|\nabla\mathcal{S}_{\mu}^{\nu}\|_{L^{2}_{\#}(\mathbb{R},L^{2}(\torus^{2}))}+   \int_{0}^{1}\int_{\torus^{2}}\frac{(\tilde{\mathcal{A}}_{\epsilon}+\nu)}{2\epsilon}
|\Delta\mathcal{S}_{\mu}^{\nu}|^{2}dxd\theta\leq\frac{\epsilon\gamma^{2}}{\nu}\Big(\int_{0}^{1}\int_{\torus^{2}}|\nabla\mathcal{S}_{\mu}^{\nu}|^{2}dxd\theta+1\Big).$$
From this last inequality, we deduce  
$$\frac{\nu}{2\epsilon}\|\Delta\mathcal{S}_{\mu}^{\nu}\|_{L^{2}_{\#}(\mathbb{R},L^{2}(\torus^{2}))}^{2}\leq\frac{\epsilon\gamma^{2}}{\nu}\big(\frac{\gamma^2}{\nu^2}+1\big),$$ then 
$$\|\Delta\mathcal{S}_{\mu}^{\nu}\|_{L^{2}_{\#}(\mathbb{R},L^{2}(\torus^{2}))}^{2}\leq\frac{2\epsilon^2\gamma^{2}}{\nu^{2}}\big(\frac{\gamma^2}{\nu^2}+1\big),$$
which gives (\ref{+02}).\\
As $\|\nabla\mathcal{S}_{\mu}^{\nu}\|_{L^{2}_{\#}(\mathbb{R},L^{2}(\torus^{2}))}$ is bounded by $\frac{\gamma}{\nu}$ (see (\ref{+01})), we can deduce that there exists a $\theta_{0}\in[0,1]$ such that 
\begin{equation}\|\nabla\mathcal{S}_{\mu}^{\nu}(\theta_{0},\cdot )\|_{2}\leq\frac{\gamma}{\nu}.\end{equation}
From (\ref{L15}) we have 
\begin{equation}\label{L150}\frac{d\|\nabla\mathcal{S}_{\mu}^{\nu}\|_{2}^{2}}{d\theta}\leq \frac{2\epsilon\gamma^{2}}{\nu}\big(\int_{\torus^{2}}|\nabla\mathcal{S}_{\mu}^{\nu}|^{2}dx+1\big).
\end{equation}
Integrating (\ref{L150}) from $\theta_{0}$ to an other $\theta_{1}\in[0,1]$ gives

\begin{multline}\label{3.55}
    \|\nabla\mathcal{S}_{\mu}^{\nu}(\theta_{1},\cdot)\|_{2}^{2}-\|\nabla\mathcal{S}_{\mu}^{\nu}(\theta_{0},\cdot)\|_{2}^{2}
\leq\frac{2\epsilon\gamma^{2}}{\nu}\int_{\theta_{0}}^{\theta_{1}}
    \left(\int_{\torus^{2}}|\nabla\mathcal{S}_{\mu}^{\nu}|^{2}dx+1\right)d\theta
\\
\leq \frac{2\epsilon\gamma^{2}}{\nu}\left(\|\nabla\mathcal{S}_{\mu}^{\nu}\|^2_{L^{2}_{\#}(\mathbb{R},L^{2}(\torus^{2}))} +1\right),~~
\end{multline}
giving the sought bound on $\|\nabla\mathcal{S}_{\mu}^{\nu}(\theta_{1},\cdot)\|_{L^{\infty}_{\#}(\mathbb{R},L^{2}(\torus^{2}))}$ for any $\theta_{1}$ or, in other words (\ref{3.282}). 
\\
Using Fourier expansion argument, because of (\ref{3.280}), we have
\begin{equation}\label{3.60}
    \| \mathcal{S}_{\mu}^{\nu}(\theta,\cdot)\|_{2}^{2}\leq \| \nabla\mathcal{S}_{\mu}^{\nu}(\theta,\cdot)\|_{2}^{2}\leq \frac{\gamma^2}{\nu^2}+2\frac{\epsilon\gamma^2}{\nu}(\frac{\gamma^2}{\nu^2}+1),
\end{equation}
and then (\ref{3.283}).\\
We have that $\frac{\partial\mathcal{S}_{\mu}^{\nu}}{\partial t}$ is solution to 
\begin{equation}\label{eqa2}
\mu\frac{\partial\mathcal{S}_{\mu}^{\nu}}{\partial t}+\frac{\partial\Big(\frac{\partial\mathcal{S}_{\mu}^{\nu}}{\partial t}\Big)}{\partial\theta}-\frac{1}{\epsilon}\nabla\cdot\Big(\big(\tilde{\mathcal{A}}_{\epsilon}+\nu\big)\nabla\big(\frac{\partial\mathcal{S}_{\mu}^{\nu}}{\partial t}\big)\Big)=\frac{1}{\epsilon}\nabla\cdot\Big(\frac{\partial\tilde{\mathcal{C}}_{\epsilon}}{\partial t}\Big)+\frac{1}{\epsilon}\nabla\cdot\Big(\frac{\partial\tilde{\mathcal{A}}_{\epsilon}}{\partial  t}\nabla\mathcal{S}_{\mu}^{\nu}\Big),\end{equation}
from which we deduce
\begin{equation}\label{eqat2}\mu\Big\|\frac{\partial\mathcal{S}_{\mu}^{\nu}}{\partial t}\Big\|_{2}^{2}+\frac{1}{2}\frac{d\Big\|\frac{\partial\mathcal{S}_{\mu}^{\nu}}{\partial t}\Big\|_{2}^{2}}{d\theta}+\frac{1}{\epsilon}\int_{\torus^{2}}\big(\tilde{\mathcal{A}}_{\epsilon}+\nu\big)\Big|\nabla\big(\frac{\partial\mathcal{S}_{\mu}^{\nu}}{\partial t}\big)\Big|^{2}dx=-\frac{1}{\epsilon}\int_{\torus^{2}}\check{\mathcal{C}}_{\epsilon}\cdot\nabla\big(\frac{\partial\mathcal{S}_{\mu}^{\nu}}{\partial t}\big)dx,\end{equation}
where 
\begin{equation}\check{\mathcal{C}}_{\epsilon}=\frac{\partial\tilde{\mathcal{C}}_{\epsilon}}{\partial t}+\frac{\partial\tilde{\mathcal{A}}_{\epsilon}}{\partial  t}\nabla\mathcal{S}_{\mu}^{\nu},\quad
\ds\nabla\cdot\check{\mathcal{C}}_{\epsilon}=\nabla\cdot\Big(\frac{\partial\tilde{\mathcal{C}}_{\epsilon}}{\partial t}+\frac{\partial\tilde{\mathcal{A}}_{\epsilon}}{\partial  t}\nabla\mathcal{S}_{\mu}^{\nu}\Big).
\end{equation}
From (\ref{Hy1}), (\ref{+01}) and (\ref{+02}), we have
\begin{equation}
\Big\|\check{\mathcal{C}}_{\epsilon}\Big\|^{2}_{L_{\#}^{2}(\mathbb{R},L^{2}(\torus^{2}))}\leq\epsilon^2\gamma(1+\frac{\gamma}{\nu}),\,\,\Big\|\nabla\cdot\check{\mathcal{C}}_{\epsilon}\Big\|_{L^{2}_{\#}(\mathbb{R},L^{2}(\torus^{2}))}\leq\epsilon^2\gamma\Big(1+\frac{\gamma}{\nu}+\epsilon\sqrt{\epsilon}\frac{\gamma}{\nu}\sqrt{\frac{\gamma}{\nu}+1}\Big).
\end{equation}
Integrating (\ref{eqat2}) from $0$ to  $1$ with respect to the variable  $\theta,$ we obtain  
$$\frac{\nu}{\epsilon}\Big\|\nabla\frac{\partial\mathcal{S}_{\mu}^{\nu}}{\partial t}\Big\|_{L^{2}_{\#}(\mathbb{R},L^{2}(\torus^{2}))}^{2}\leq\epsilon^2\gamma(1+\frac{\gamma}{\nu})\Big\|\nabla\frac{\partial\mathcal{S}_{\mu}^{\nu}}{\partial t}\Big\|_{L^{2}_{\#}(\mathbb{R},L^{2}(\torus^{2}))},$$ 
then
$$\Big\|\nabla\frac{\partial\mathcal{S}_{\mu}^{\nu}}{\partial t}\Big\|_{L^{2}_{\#}(\mathbb{R},L^{2}(\torus^{2}))}\leq\epsilon^3\frac{\gamma}{\nu}(1+\frac{\gamma}{\nu}).$$
Using the Fourier expansion of  $\mathcal{S}_{\mu}^{\nu},$ we have for a given $\theta_{0}$ 
$$\Big\|\nabla\frac{\partial\mathcal{S}_{\mu}^{\nu}}{\partial t}(\theta_{0},\cdot)\Big\|_{2}\leq\epsilon^3\frac{\gamma}{\nu}(1+\frac{\gamma}{\nu}).$$
Thus,  as previously, we get
$$\Big\|\nabla\frac{\partial\mathcal{S}_{\mu}^{\nu}}{\partial t}\Big\|_{L^{\infty}_{\#}(\mathbb{R},L^{2}(\torus^{2}))}\leq\epsilon^3\frac{\gamma}{\nu}(1+\frac{\gamma}{\nu}),\,\,\,\,\Big\|\frac{\partial\mathcal{S}_{\mu}^{\nu}}{\partial t}\Big\|_{L^{\infty}_{\#}(\mathbb{R},L^{2}(\torus^{2}))}\leq\epsilon^3\frac{\gamma}{\nu}(1+\frac{\gamma}{\nu}).$$
\end{proof}
Since the estimates of theorem \ref{theoexes12} do not depend on $\mu,$ making the process $\mu\rightarrow0$ allows us to deduce the following theorem.

\begin{theorem}\label{theoexes02}
Under assumptions (\ref{Hy1}),(\ref{Hy2}) and (\ref{Hy3}), for any $\nu>0,$ there exists a unique $\mathcal{S}^{\nu}=\mathcal{S}^{\nu}(t,\theta,x)\in L^{2}(\mathbb{R}\times\torus^{2}),$ periodic of period 1 with respect to $\theta$ solution to (\ref{L6}) and submitted to the constraint
\begin{equation}\label{3.710}
    \sup_{\theta\in\mathbb{R}}\Big|\int_{\torus^{2}}\mathcal{S}^{\nu}(\theta,x)dx\Big|=0.
\end{equation}
Moreover, the following estimates are satisfied
\begin{equation}\label{3.711}\begin{array}{ccc}
    \ds\Big\|\frac{\partial\mathcal{S}^{\nu}}{\partial\theta}\Big\|_{L^{2}_{\#}(\mathbb{R},L^{2}(\torus^{2}))}\leq\frac{\gamma}{\sqrt{\epsilon\nu}}\sqrt{(\frac{\gamma}{2\nu}+1)},\qquad
    \|\nabla\mathcal{S}^{\nu}\|_{L^{\infty}_{\#}(\mathbb{R},L^{2}(\torus^{2}))}\leq\sqrt{\frac{\gamma^2}{\nu^2}+\frac{2\epsilon\gamma^2}{\nu}(\frac{\gamma^2}{\nu^2}+1)},
   \end{array}\end{equation}
\begin{equation}\label{3.712}\ds \|\mathcal{S}^{\nu}\|_{L^{\infty}_{\#}(\mathbb{R},L^{2}(\torus^{2}))}\leq\sqrt{\frac{\gamma^2}{\nu^2}+\frac{2\epsilon\gamma^2}{\nu}(\frac{\gamma^2}{\nu^2}+1)},
\qquad
    \Big\|\frac{\partial\mathcal{S}^{\nu}}{\partial t}\Big\|_{L^{\infty}_{\#}(\mathbb{R},L^{2}(\torus^{2}))}\leq\epsilon^3\frac{\gamma}{\nu}(1+\frac{\gamma}{\nu}).\end{equation}
\end{theorem}

~
\,\,\,\,\begin{proof}.\,\,(of Theorem \ref{theoexes02}). As estimates of Theorem \ref{theoexes12} do not depend on $\mu,$
to proof existence of $\mathcal{S}^{\nu},$ it suffices to make $\mu$ tend to $0$ in (\ref{L7}).
%
Uniqueness is insured by (\ref{3.710}), once noticed that, if $\mathcal{S}^{\nu}$ and $\widetilde{\mathcal{S}}^{\nu}$ are two solutions of (\ref{L6}), with constraint (\ref{3.710}), $\mathcal{S}^{\nu}-\widetilde{\mathcal{S}}^{\nu}$ is solution to
\begin{equation}\label{3.713}
\frac{\partial(\mathcal{S}^{\nu}-\widetilde{\mathcal{S}}^{\nu})}{\partial\theta}-\frac{1}{\epsilon}\nabla\cdot((\widetilde{\mathcal{A}}_{\epsilon}+\nu)
\nabla(\mathcal{S}^{\nu}-\widetilde{\mathcal{S}}^{\nu}))=0,
    \end{equation}
from which we can deduce that
\begin{equation}\label{3.714}
    \nu\|\nabla(\mathcal{S}^{\nu}-\widetilde{\mathcal{S}}^{\nu})\|^{2}_{L^{2}_{\#}(\mathbb{R},L^{2}(\torus^{2}))}=0,
\end{equation}
and  because of (\ref{3.710}), and its consequence:
\begin{equation}\label{3.715}
    \|\mathcal{S}^{\nu}-\widetilde{\mathcal{S}}^{\nu}\|_{L^{2}_{\#}(\mathbb{R},L^{2}(\torus^{2}))}\leq
    \|\nabla(\mathcal{S}^{\nu}-\widetilde{\mathcal{S}}^{\nu})\|_{L^{2}_{\#}(\mathbb{R},L^{2}(\torus^{2}))} ,
\end{equation}
that
\begin{equation}\label{3.716}
    \widetilde{\mathcal{S}}^{\nu}=\mathcal{S}^{\nu}.
\end{equation}
\end{proof}
Now we get estimates on $\mathcal{S}^\nu$ which do not depend on $\nu.$
\begin{theorem}\label{th2.3}
Under the assumptions (\ref{Hy1}),(\ref{Hy2}) and (\ref{Hy3}), the solution $\mathcal{S}^{\nu},$ of (\ref{L6}) given by theorem \ref{theoexes02} satisfies the following properties
\begin{equation}\label{gr2}
\Big\|\sqrt{\widetilde{\mathcal{A}}_{\epsilon}}|\nabla\mathcal{S}^{\nu}|\Big\|_{L^{2}_{\#}(\mathbb{R},L^{2}(\torus^{2}))}\,\,\leq \gamma,
\end{equation}
\begin{equation}\label{gr3}
\Big(\int_{\theta_{\alpha}}^{\theta_{\omega}}\int_{\torus^{2}}|\nabla\mathcal{S}^{\nu}|^{2}dxd\theta \Big)^{1/2}\,\,\leq\,\,\frac{\gamma}{
\sqrt{\widetilde{G}_{thr}}},
\end{equation}
\begin{equation}\label{gr3bis}
\Big\|\nabla\mathcal{S}^{\nu}(\theta_0,\cdot)\Big\|_{2}\,\,\leq\,\,\frac{\gamma}{
\sqrt{\widetilde{G}_{thr}}},\,\,\text{for a given}\,\,\theta_0\in\,\,[\theta_\alpha,\theta_\omega],
\end{equation}
\begin{equation}\label{XXC}\|\mathcal{S}^{\nu}\|^{2}_{L_{\#}^\infty(\mathbb{R},L^2(\torus^2))}\leq\frac{\gamma}{\sqrt{\tilde{G}_{thr}}}+ 2\epsilon\gamma^3,
\end{equation}
\begin{equation}\label{3.1013}
    \left\|\frac{\partial\mathcal{S}^{\nu}}{\partial t}\right\|^{2}_{L_{\#}^\infty(\mathbb{R},L^2(\torus^2))}\leq\,\epsilon\Big(\frac{\gamma+\epsilon\gamma^{3}}{\sqrt{\widetilde G_{thr}}}+(\gamma^2+\epsilon^2\gamma^4)\Big).
\end{equation}

\end{theorem}
~
\,\,\,\,\begin{proof}.\,(of Theorem \ref{th2.3}).\,\,
Multiplying (\ref{L6}) by $\mathcal{S}^{\nu}$ and integrating over $\torus^{2}$ yields 
\begin{equation}\label{gr1}
\frac{1}{2}\frac{d}{d\theta}\int_{\torus^{2}}|\mathcal{S}^{\nu}|^{2}dx+\frac{1}{\epsilon}\int_{\torus^{2}}(\widetilde{\mathcal{A}}_{\epsilon}+\nu)|\nabla\mathcal{S}^{\nu}|^{2}dx=-\frac{1}{\epsilon}\int_{\torus^{2}}\widetilde{\mathcal{C}}_{\epsilon}\cdot\nabla\mathcal{S}^{\nu}dx.\end{equation}
Integrating  (\ref{gr1}) in $\theta$ over $[0,1]$  gives 
\begin{equation}
\frac{1}{\epsilon}\int_{0}^{1}\int_{\torus^{2}}(\widetilde{\mathcal{A}}_{\epsilon}+\nu)|\nabla\mathcal{S}^{\nu}|^{2}dx\leq\frac{\gamma}{\epsilon}\int_{0}^{1}\int_{\torus^{2}}\sqrt{\widetilde{\mathcal{A}}_{\epsilon}}|\nabla\mathcal{S}^{\nu}|dx,
\end{equation}
then we obtain (\ref{gr2}).\\
Assuming (\ref{Hy2}),  we have
\begin{equation}
\sqrt{\widetilde{G}_{thr}}\Big(\int_{\theta_{\alpha}}^{\theta_{\omega}}\int_{\torus^{2}}|\nabla\mathcal{S}^{\nu}|^{2}dxd\theta \Big)^{1/2}\leq\Big( \int_{\theta_{\alpha}}^{\theta_{\omega}}\int_{\torus^{2}}\widetilde{\mathcal{A}}_{\epsilon}|\nabla\mathcal{S}^{\nu}|^{2}dxd\theta\Big)^{\frac{1}{2}}\leq\Big\|\sqrt{\widetilde{\mathcal{A}}_{\epsilon}}|\nabla\mathcal{S}^{\nu}|\Big\|_{L^{2}_{\#}(\mathbb{R},L^{2}(\torus^{2}))}.
\end{equation}
From (\ref{gr2}) and this last inequality we get (\ref{gr3}).
Then, there exists  a $\theta_{0}\in[\theta_{\alpha},\theta_{\omega}]$ such that $\mathcal{S}^{\nu}$ satisfies (\ref{gr3bis}).\\
Using the Fourier expansion of $\mathcal{S}^{\nu}$ and the relation (\ref{3.710})  we get 
\begin{equation}\label{gr6}\Big\|\mathcal{S}^{\nu}(\theta_{0},\cdot)\Big\|_{2}\leq \Big\|\nabla\mathcal{S}^{\nu}(\theta_{0},\cdot)\Big\|_{2}\leq\frac{\gamma}{\sqrt{\tilde{G}_{thr}}}.\end{equation}
%
Multiplying (\ref{L6}) by $\mathcal{S}^{\nu},$ integrating over $\torus^{2}$ we obtain

$$\frac{1}{2}\frac{d\|\mathcal{S}^{\nu}(\theta,\cdot)\|^{2}_{2}}{d\theta}+\frac{1}{\epsilon}\int_{\torus^{2}}(\widetilde{\mathcal{A}}_{\epsilon}+\nu)|\nabla\mathcal{S}^{\nu}(\theta,\cdot)|^{2}dx=\frac{1}{\epsilon}\int_{\torus^{2}}|\nabla\cdot\widetilde{\mathcal{C}}_{\epsilon}\mathcal{S}^{\nu}(\theta,\cdot)|dx$$
Applying formula (\ref{eqa})
with $V=\frac{|\nabla\cdot\widetilde{\mathcal{C}}_\epsilon|}{\epsilon}$ and $U=|\mathcal{S}^\nu|,$ we get
$$\frac{1}{2}\frac{d\|\mathcal{S}^{\nu}(\theta,\cdot)\|^{2}_{2}}{d\theta}+\frac{1}{\epsilon}\int_{\torus^{2}}(\widetilde{\mathcal{A}}_{\epsilon}+\nu)|\nabla\mathcal{S}^{\nu}(\theta,\cdot)|^{2}dx\leq\int_{\torus^{2}}\Big[\frac{(\widetilde{\mathcal{A}}_{\epsilon}+\nu)}{4\epsilon}|\mathcal{S}^{\nu}(\theta,\cdot)|^2+\frac{1}{\epsilon(\widetilde{\mathcal{A}}_{\epsilon}+\nu)}|\nabla\cdot\widetilde{\mathcal{C}}_{\epsilon}|^2\Big]\,dx,$$
which gives
\begin{equation}\label{SSX}\frac{1}{2}\frac{d\|\mathcal{S}^{\nu}(\theta,\cdot)\|^{2}_{2}}{d\theta}+\frac{1}{\epsilon}\int_{\torus^{2}}(\widetilde{\mathcal{A}}_{\epsilon}+\nu)\Big(|\nabla\mathcal{S}^{\nu}(\theta,\cdot)|^{2}-\frac{|\mathcal{S}^{\nu}(\theta,\cdot)|^2}{4}\Big)dx\leq\int_{\torus^{2}}\frac{1}{\epsilon(\widetilde{\mathcal{A}}_{\epsilon}+\nu)}|\nabla\cdot\widetilde{\mathcal{C}}_{\epsilon}|^2\,dx.
\end{equation}
Using Fourier expansion of $\mathcal{S}^{\nu}(\theta,\cdot),$ one can prove that the second term of the left hand side of (\ref{SSX}) is positive, then we have
\begin{equation}\label{fra}
\frac{d\|\mathcal{S}^{\nu}(\theta,\cdot)\|^{2}_{2}}{d\theta}\leq2\int_{\torus^{2}}\frac{1}{\epsilon(\widetilde{\mathcal{A}}_{\epsilon}+\nu)}|\nabla\cdot\widetilde{\mathcal{C}}_{\epsilon}|^2\,dx
.
\end{equation}
Using (\ref{Hy1}), (\ref{Hy3}) and integrating (\ref{fra}) from $\theta_0$ to $\theta\in[0,1],$ we obtain
\begin{equation}\|\mathcal{S}^{\nu}(\theta,\cdot)\|_{2}^{2}\leq\|\mathcal{S}^{\nu}(\theta_{0},\cdot)\|_{2}^{2}+2\epsilon\gamma^3,\end{equation} 
then  inequality (\ref{XXC}) is satisfied.\\
Using inequality (\ref{gr2}) and from hypothesis (\ref{Hy3}) we get
\begin{equation} \Big\|\frac{\partial (\nabla\widetilde{\mathcal{A}}_{\epsilon})}{\partial t}\nabla\mathcal{S}^{\nu}\Big\|_{L_{\#}^{2}(\mathbb{R},L^{2}(\torus^{2}))}\,\,\leq\,\,\epsilon^2\gamma\Big\|\sqrt{\widetilde{\mathcal{A}}_{\epsilon}}\nabla\mathcal{S}^{\nu}\Big\|_{L_{\#}^{2}(\mathbb{R},L^{2}(\torus^{2}))}\leq\,\,\epsilon^2\gamma^{2}.\end{equation}
Multiplying (\ref{L6}) by $-\Delta\mathcal{S}^{\nu}$ and integrating in $x\in \torus^{2}$ we get
\begin{equation}\label{gr8}
\frac{1}{2}\frac{d}{d\theta}\|\nabla\mathcal{S}^{\nu}\|_{2}^{2}+\frac{1}{\epsilon}\int_{\torus^{2}}\Big(\widetilde{\mathcal{A}}_{\epsilon}+\nu\Big)|\Delta\mathcal{S}^{\nu}|^{2}dx
+\frac{1}{\epsilon}\int_{\torus^{2}}\nabla\widetilde{\mathcal{A}}_{\epsilon}\cdot\nabla\mathcal{S}^{\nu}\Delta\mathcal{S}^{\nu}dx=-\frac{1}{\epsilon}\int_{\torus^{2}}\nabla\cdot\widetilde{\mathcal{C}}_{\epsilon}\cdot\Delta\mathcal{S}^{\nu}dx.
\end{equation}
Using (\ref{eqa}) with $U=|\Delta\mathcal{S}^{\nu}|$ and $V= \frac{\nabla\widetilde{\mathcal{A}}_{\epsilon}\cdot\nabla\mathcal{S}^{\nu}}{\epsilon}$ and with $U=|\Delta\mathcal{S}^{\nu}|$ and $V=\frac{\nabla\cdot\widetilde{\mathcal{C}}_{\epsilon}}{\epsilon},$
the equality (\ref{gr8}) becomes
\begin{equation}
\frac{1}{2}\frac{d}{d\theta}\|\nabla\mathcal{S}^{\nu}\|_{2}^{2}+\frac{1}{2\epsilon}\int_{\torus^{2}}(\widetilde{\mathcal{A}}_{\epsilon}+\nu)|\Delta\mathcal{S}^{\nu}|^{2}dx\leq \,\,\frac{1}{\epsilon}\int_{\torus^{2}}\Big[\frac{|\nabla\widetilde{\mathcal{A}}_{\epsilon}|^{2}}{\widetilde{\mathcal{A}}_{\epsilon}+\nu}|\nabla\mathcal{S}^{\nu}|^{2}+\frac{|\nabla\widetilde{\mathcal{C}}_{\epsilon}|^{2}}{\widetilde{\mathcal{A}}_{\epsilon}+\nu}\Big]dx,
\end{equation}
which, integrating from $0$ to $1$ yields
\begin{equation}
\int_{0}^{1}\int_{\torus^{2}}\widetilde{\mathcal{A}}_{\epsilon}|\Delta\mathcal{S}^{\nu}|^{2}dxd\theta\leq\,\,2\epsilon\gamma^2\Big(\int_{0}^{1}\int_{\torus^{2}}
|\widetilde{\mathcal{A}}_{\epsilon}||\nabla\mathcal{S}^{\nu}|^{2}dxd\theta+\gamma\Big) \,\,\leq\,\,2\epsilon\gamma^2(\gamma^{2}+\gamma).
\end{equation}
As \begin{equation}
\Big|\frac{\partial\widetilde{\mathcal{A}}_{\epsilon}}{\partial t}\Big|\,\,\leq\,\,\epsilon^2\gamma|\widetilde{\mathcal{A}}_{\epsilon}|,
\end{equation} we obtain
\begin{equation}
\Big\|\sqrt{\frac{\partial\widetilde{\mathcal{A}}_{\epsilon}}{\partial t}} \Delta\mathcal{S}^{\nu}\Big\|_{L_{\#}^{2}(\mathbb{R},L^{2}(\torus^{2}))}
\,\,\leq\,\,\epsilon\sqrt{2\epsilon}\gamma^2\sqrt{1+\gamma}.\end{equation}
Now we set out the equation to which $\frac{\partial \mathcal{S}^{\nu}}{\partial t}$ is solution. We have $$\frac{\partial }{\partial\theta}\Big(\frac{\partial \mathcal{S}^{\nu}}{\partial t} \Big)=\frac{\partial }{\partial t}\Big(\frac{\partial \mathcal{S}^{\nu}}{\partial\theta} \Big)=\frac{1}{\epsilon}\Big(\nabla\cdot\frac{\partial \widetilde{\mathcal{A}}_{\epsilon}}{\partial t}\nabla\mathcal{S}^{\nu}+(\widetilde{\mathcal{A}}_{\epsilon}+\nu)\nabla\frac{\partial \mathcal{S}^{\nu}}{\partial t}\Big)+\frac{1}{\epsilon}\nabla\cdot\frac{\partial \widetilde{\mathcal{A}}_{\epsilon}}{\partial t},$$
then $\frac{\partial \mathcal{S}^{\nu}}{\partial t}$ is solution to
\begin{equation}\label{gr9}\frac{\partial }{\partial\theta}\Big(\frac{\partial \mathcal{S}^{\nu}}{\partial t}\Big)-\frac{1}{\epsilon}\nabla\cdot\Big((\widetilde{\mathcal{A}}_{\epsilon}+\nu)\nabla\big(\frac{\partial\mathcal{S}^{\nu}}{\partial t}\big)\Big)=\frac{1}{\epsilon}\nabla\cdot\Big(\frac{\partial \widetilde{\mathcal{A}}_{\epsilon}}{\partial t}\nabla\mathcal{S}^{\nu}\Big)+\frac{1}{\epsilon}\nabla\cdot\frac{\partial \widetilde{\mathcal{C}}_{\epsilon}}{\partial t}.
\end{equation}
Multiplying (\ref{gr9}) by $\frac{\partial \mathcal{S}^{\nu}}{\partial t}$ and integrating in $x\in\torus^{2},$ we get
\begin{equation}\label{gr10}
\frac{1}{2}\frac{d}{d\theta}\Big\|\frac{\partial \mathcal{S}^{\nu}}{\partial t}\Big\|_{2}^{2}+\frac{1}{\epsilon}\int_{\torus^{2}}(\widetilde{\mathcal{A}}_{\epsilon}+\nu)\Big|\nabla\frac{\partial\mathcal{S}^{\nu}}{\partial t}\Big|^{2}dx\leq \,\,\frac{1}{\epsilon}\int_{\torus^{2}}\Big|\frac{\partial\widetilde{\mathcal{A}}_{\epsilon}}{\partial t}\Big||\nabla\mathcal{S}^{\nu}|\Big|\nabla\frac{\partial \mathcal{S}^{\nu}}{\partial t}\Big|dx+\frac{1}{\epsilon}\int_{\torus^{2}}\Big|\frac{\partial\widetilde{\mathcal{C}}_{\epsilon}}{\partial t}\Big|\Big|\nabla\frac{\partial \mathcal{S}^{\nu}}{\partial t}\Big|dx.
\end{equation}
Using the fact that $\Big|\frac{\partial\widetilde{\mathcal{C}}_{\epsilon}}{\partial t}\Big|^{2}\,\,\leq\,\,\epsilon^2\gamma^{2}\big|\widetilde{\mathcal{A}}_{\epsilon}\big|, $ the second term of the right hand side of (\ref{gr10}) satisfies
 \begin{eqnarray}\label{3.95a}
 {    \int_{\torus^{2}}\bigg|\frac{\partial\widetilde{\mathcal{C}}_{\epsilon}}{\partial t}\bigg|\;\bigg|\nabla\frac{\partial\mathcal{S}^{\nu}}{\partial t}\bigg|dx\leq\epsilon\gamma\int_{\torus^{2}}\sqrt{\widetilde{\mathcal{A}}_{\epsilon}}\bigg|\nabla\frac{\partial\mathcal{S}^{\nu}}{\partial t}\bigg|dx {} }
{}    \leq\epsilon\gamma\bigg\|\sqrt{\widetilde{\mathcal{A}}_{\epsilon}}\;\left|\nabla\frac{\partial\mathcal{S}^{\nu}}{\partial t}\right|\bigg\|_{2}.
\end{eqnarray}
In the same way, using (\ref{Hy3}) we deduce the following estimate for the first term of the right hand side of (\ref{gr10})
 \begin{eqnarray}\label{3.96a}{
 \int_{\torus^{2}}\left|\frac{{\partial \widetilde{\mathcal{A}}}_{\epsilon}}{\partial t}\right| \,|\nabla\mathcal{S}^{\nu}| \left|\nabla\frac{\partial\mathcal{S}^{\nu}}{\partial t}\right|dx\leq
 \left\|\sqrt{\bigg|\frac{\partial\widetilde{\mathcal{A}}_{\epsilon}}{\partial t}}\bigg|\, |\nabla\mathcal{S}^{\nu}|\right\|_{2}\,\left\|\sqrt{\bigg|\frac{\partial\widetilde{\mathcal{A}}_{\epsilon}}{\partial t}}\bigg| \Big|\nabla\frac{\partial\mathcal{S}^{\nu}}{\partial t}\Big|\right\|_{2}{} }~~~~~~~~
 \nonumber\\
{}
\leq \epsilon^2\gamma^{2}\left\|\sqrt{\widetilde{{\mathcal{A}}}_{\epsilon}} \left|\nabla\mathcal{S}^{\nu}\right|\right\|_{2}\,
\left\|\sqrt{\widetilde{{\mathcal{A}}}_{\epsilon}}\left| \nabla\frac{\partial\mathcal{S}^{\nu}}{\partial t}\right|\right\|_{2}.
\end{eqnarray}
Using inequalities  (\ref{3.95a}), (\ref{3.96a}) and (\ref{gr2}) and integrating (\ref{gr10}) in $\theta$ over $[0,1]$, we have
\begin{eqnarray}\label{3.97a}{
\Bigg\|\sqrt{(\widetilde{{\mathcal{A}}}_{\epsilon}+\nu)}\,\bigg|\nabla\frac{\partial\mathcal{S}^{\nu}}{\partial t}\bigg|\Bigg\|_{L^{2}_
{\#}(\mathbb{R},L^{2}(\torus^{2}))}^{2}
\,\,\leq\,\,\epsilon\gamma\Bigg\|\sqrt{\widetilde{{\mathcal{A}}}_{\epsilon}}
\bigg|\nabla\frac{\partial\mathcal{S}^{\nu}}{\partial t}\bigg|\Bigg\|_{L^{2}_{\#}(\mathbb{R},L^{2}(\torus^{2}))} {} }~~~~~~~~~~~~~~~
\nonumber\\
+\epsilon^2\gamma^{3}\Bigg\|\sqrt{\widetilde{{\mathcal{A}}}_{\epsilon}}\bigg|\nabla\frac{\partial\mathcal{S}^{\nu}}{\partial t}\bigg |\Bigg\|_{L^{2}_{\#}(\mathbb{R},L^{2}(\torus^{2}))}.
\end{eqnarray}
From this last inequality, we deduce
\begin{equation}\label{3.98a}
    \left\|\sqrt{\widetilde{{\mathcal{A}}}_{\epsilon}}\left|\nabla\frac{\partial\mathcal{S}^{\nu}}{\partial t}\right|\right\|_{L^{2}_{\#}(\mathbb{R},L^{2}(\torus^{2}))}\leq \epsilon(\gamma+\epsilon\gamma^{3}),
\end{equation}
and then
\begin{equation}\label{3.909}
    \int_{\theta_{\alpha}}^{\theta_{\omega}}\left\|\nabla\frac{\partial\mathcal{S}^{\nu}}{\partial t}\right\|_{2}d\theta\,\,\leq\,\,\epsilon\frac{\gamma+\epsilon\gamma^{3}}{\sqrt{\widetilde G_{thr}}}.
\end{equation}
From (\ref{3.909}), we deduce that there exists a $\theta_{0}\in[\theta_{\alpha},\theta_{\omega}]$ such that
\begin{equation}\label{3.100}
    \left\|\nabla\frac{\partial\mathcal{S}^{\nu}}{\partial t}(\theta_{0},\cdot)\right\|_{2}\leq\epsilon\frac{\gamma+\epsilon\gamma^{3}}{\sqrt{\widetilde G_{thr}}},
\end{equation}
and, since the mean value of $\frac{\partial\mathcal{S}^{\nu}}{\partial t}(\theta_{0},\cdot)$
is zero,
\begin{equation}\label{3.101}
    \left\|\frac{\partial\mathcal{S}^{\nu}}{\partial t}(\theta_{0},\cdot)\right\|_{2}\,\,\leq\,\,\epsilon\frac{\gamma+\epsilon\gamma^{3}}{\sqrt{\widetilde G_{thr}}}.
\end{equation}
To end the proof of the theorem it remains to show that $\frac{\partial \mathcal{S}^{\nu}}{\partial t}$ is bounded independently of $\nu$ in $L^{\infty}_{\#}(\mathbb{R},L^{2}(\torus^{2})).$ For this we will estimate the right hand side of (\ref{gr10}) by applying formula (\ref{eqa}) with $V=\frac{1}{\epsilon}|\frac{\partial\widetilde{\mathcal{C}}_{\epsilon}}{\partial t}|$ and $U=|\nabla\frac{\partial\mathcal{S}^{\nu}}{\partial t}|$ to treat the second term of the right hand side of (\ref{gr10}) and with  $V=\frac{1}{\epsilon}|\frac{\partial\widetilde{\mathcal{A}}_{\epsilon}}{\partial t}||\nabla\mathcal{S}^{\nu}|$ and $U=|\nabla\frac{\partial\mathcal{S}^{\nu}}{\partial t}|$
to treat the first. It gives:
\begin{multline}\label{3.102}{\frac{1}{2}\frac{\ds\partial\bigg(\Big\|\frac{\ds\partial\mathcal{S}^{\nu}}{\partial t}\Big\|_{2}^{2}\bigg)}{\partial\theta}+\int_{\torus^{2}}\frac{\widetilde{\mathcal{A}}_{\epsilon}+\nu}{2\epsilon}\left|\nabla\frac{\partial\mathcal{S}^{\nu}}{\partial t}\right|^{2}dx{}}
\\
{}\leq\,\,\int_{\torus^{2}}\frac{\ds\Big|\frac{\partial\widetilde{\mathcal{C}}_{\epsilon}}{\partial t}\Big|^2}{\epsilon(\widetilde{\mathcal{A}}_{\epsilon}+\nu)}dx+\int_{\torus^{2}}\frac{\ds\Big|\frac{\partial\widetilde{\mathcal{A}}_{\epsilon}}{\partial t}\Big|^{2}|\nabla\mathcal{S}^{\nu}|^{2}}{\epsilon(\widetilde{\mathcal{A}}_{\epsilon}+\nu)}dx
\leq \epsilon\gamma^2+\epsilon^3\gamma^2\int_{\torus^{2}}\widetilde{\mathcal{A}}_{\epsilon}|\nabla\mathcal{S}^{\nu}|^{2}dx,
\end{multline}
where we used hypothesis (\ref{Hy3}) to get the last inequality. Integrating this last formula in $\theta$ over $[\theta_{0},\sigma]$ for any $\sigma>\theta_{0},$ we obtain, always remembering (\ref{gr2}),
\begin{equation}\label{3.1013a}
    \left\|\frac{\partial\mathcal{S}^{\nu}}{\partial t}(\sigma,\cdot)\right\|_{2}^{2}\leq \left\|\frac{\partial\mathcal{S}^{\nu}}{\partial t}(\theta_{0},\cdot)\right\|_{2}^{2}+\epsilon(\gamma^2+\epsilon^2\gamma^4).
\end{equation}
From inequality (\ref{3.1013a}) we obtain directly the   inequality of (\ref{3.1013}), using the periodicity of $\mathcal{S}^{\nu}.$ 
\end{proof}

Estimates (\ref{XXC}) and (\ref{3.1013}) given in theorem \ref{th2.3} do not depend on $\nu.$ Making $\nu\rightarrow 0,$ allows us to deduce that, up to a subsequence $\,\,\mathcal{S}^{\nu}\longrightarrow S\in L_\#^\infty(\mathbb{R},L^2(\torus^2))\,\,\textrm{weak}-*.$ Concerning the limit $\mathcal{S}$ we have the following theorem. 
\begin{theorem}\label{th3.017}
Under assumptions (\ref{Hy1}), (\ref{Hy2}), (\ref{Hy3}), 
there exists a unique function $\mathcal{S}=\mathcal{S}(t,\theta,x)\in L_{\#}^{\infty}(\mathbb{R},L^{2}(\torus^{2}))$, periodic of period 1 with respect to $\theta,$ solution to
\begin{equation}\label{3.104}
    \frac{\partial\mathcal{S}}{\partial \theta}-\frac{1}{\epsilon}\nabla\cdot(\widetilde{\mathcal{A}}_{\epsilon}(t,\cdot,\cdot)\nabla\mathcal{S})
    =\frac{1}{\epsilon}\nabla\cdot\widetilde{\mathcal{C}}_{\epsilon}(t,\cdot,\cdot),
\end{equation}
and satisfying, for any $t,\theta\in\mathbb{R}^{+}\times\mathbb{R}$
\begin{equation}\label{3.104.1}
\int_{\torus^{2}}\mathcal{S}(t,\theta,x)dx=0.
\end{equation}
Moreover it satisfies: 
\begin{equation}\label{3.105}
    \|\mathcal{S}\|_{L^{\infty}_{\#}(\mathbb{R},L^{2}(\torus^{2}))}^2\leq\frac{\gamma}{\sqrt{\widetilde{G}_{thr}}}+2\epsilon\gamma^3,
\end{equation}
\begin{equation}\label{3.106}
    \Big\|\frac{\partial\mathcal{S}}{\partial t}\Big\|^{2}_{L^{\infty}_{\#}(\mathbb{R},L^{2}(\torus^{2}))}\leq\epsilon\Big(\frac{\gamma+\epsilon\gamma^{3}}{\sqrt{\widetilde G_{thr}}}+(\gamma^2+\epsilon^2\gamma^4)\Big).
\end{equation}
\end{theorem}
\begin{proof}.\,(of Theorem\,\ref{th3.017}).\,\,
Uniqueness of $\mathcal{S}$ is not gotten via the above evoked process
${\nu}\longrightarrow 0$, but directly comes from (\ref {3.104}).
Assuming that there are two solutions ${\mathcal{S}_1}$ and ${\mathcal{S}_2}$ to (\ref {3.104}),
we easily deduce that
\begin{equation}\label{3.106.1}
    \frac{d\left(\left\|{\mathcal{S}_1}-{\mathcal{S}_2}\right\|^2_2\right)}{d\theta}
    +\frac{1}{\epsilon}\int_{\torus^{2}}\widetilde{\mathcal{A}}_{\epsilon}\left|\nabla({\mathcal{S}_1}- {\mathcal{S}_2}) \right|^2 dx   =0,
\end{equation}
which gives, because of the non-negativity of $\widetilde{\mathcal{A}}_{\epsilon}$,
\begin{equation}\label{3.106.2}
    \frac{d\left(\left\|{\mathcal{S}_1}- {\mathcal{S}_2}\right\|^2_2\right)}{d\theta}
    \leq 0.
\end{equation}
From (\ref {3.106.1}) we deduce that either
\begin{equation}\label{3.106.3}
  \widetilde{\mathcal{A}}_{\epsilon}\left|\nabla({\mathcal{S}_1}- {\mathcal{S}_2}) \right|^2 \equiv0,
\end{equation}
or, for any $\theta\in \R$,
\begin{equation}\label{3.106.4}
\left\|{\mathcal{S}_1}(\theta+1,\cdot)- {\mathcal{S}_2} (\theta+1,\cdot)\right\|^2_2
< \left\|{\mathcal{S}_1}(\theta,\cdot)- {\mathcal{S}_2} (\theta,\cdot)\right\|^2_2.
\end{equation}
As (\ref {3.106.4}) is not possible because of the periodicity of ${\mathcal{S}_1}$ and ${\mathcal{S}_2}$,
we deduce that (\ref {3.106.3}) is true.  Using this last information, we deduce, for instance
\begin{equation}\label{3.106.5}
\nabla({\mathcal{S}_1}- {\mathcal{S}_2})(\theta_\omega,\cdot) \equiv 0,
\end{equation}
yielding, because of property (\ref {3.104.1}),
\begin{equation}\label{3.106.6}
\left\|{(\mathcal{S}_1}- {\mathcal{S}_2}) (\theta_\omega,\cdot)\right\|^2_2
\leq\left\|\nabla({\mathcal{S}_1}- {\mathcal{S}_2})(\theta_\omega,\cdot)\right\|^2_2.
\end{equation}
Injecting (\ref {3.106.3}) in (\ref {3.106.1}) yields
\begin{equation}\label{3.106.7}
    \frac{d\left(\left\|{\mathcal{S}_1}- {\mathcal{S}_2}\right\|^2_2\right)}{d\theta}
    = 0,
\end{equation}
and then
\begin{equation}\label{3.106.6.111}
\left\|{(\mathcal{S}_1}- {\mathcal{S}_2}) (\theta,\cdot)\right\|^2_2=0,
\end{equation}
for any $\theta\geq\theta_\omega$ and consequently or any $\theta\in\R$. This ends the proof of theorem\,\ref{th3.017}.
\end{proof}

\noindent With this theorem on hand we can get the following result concerning $z^\epsilon$ solution of equation (\ref{L1}).

\begin{theorem}\label{th3.018}
Under properties (\ref{Hy1}), (\ref{Hy2}), (\ref{Hy3}), for any $T,$ not depending on $\epsilon,$ equation (\ref{L1}), with coefficients given by (\ref{La1}) coupled with (\ref{La3}) 
and (\ref{La4}) coupled with (\ref{La5}) 
has a unique solution $z^{\epsilon}\in L^{\infty}([0,T];L^{2}(\torus^{2})).$ This solution satisfies:
\begin{equation}\label{3.1080}
    \|z^{\epsilon}\|_{L^{\infty}([0,T],L^{2}(\torus^{2}))}\,\,\leq\widetilde{\gamma}
\end{equation}
where $\widetilde{\gamma}$ is a constant which do not depend on $\epsilon.$
\end{theorem}
\begin{proof}(of Theorem\,\ref{th01}). Theorem \,\ref{th01} is a direct consequence of theorem\,\ref{th3.018}.
\end{proof}

\noindent\begin{proof}.\,(of Theorem \ref {th3.018}). To prove uniqueness,
 we consider $z_{1}^{\epsilon}$ and $z_{2}^{\epsilon}$ two solutions of
  (\ref{L1}). Their difference is then solution to
\begin{equation}\label{3.108}\left\{\begin{array}{ccc}
    \ds \frac{\partial(z_{1}^{\epsilon}-z_{2}^{\epsilon})}{\partial t} -\frac{1}{\epsilon^{2}}\nabla\cdot\left(\widetilde{\mathcal{A}}_{\epsilon}\nabla(z_{1}^{\epsilon}-z_{2}^{\epsilon})\right)=0,\,\,
    \\
   \left( z_{1}^{\epsilon}-z_{2}^{\epsilon}\right)_{|t=0}=0,
\end{array}\right.\end{equation}
and multiplying the first equation of (\ref{3.108}) by $(z_{1}^{\epsilon}-z_{2}^{\epsilon})$ and integrating with respect to $x$ gives
\begin{equation}\label{3.109}
    \frac{d\left(\|z_{1}^{\epsilon}-z_{2}^{\epsilon}\|_{2}^{2}\right)}{dt}\,\,\leq\,\,0,
\end{equation}
yielding
\begin{equation}\label{3.110}
    \|z_{1}^{\epsilon}-z_{2}^{\epsilon}\|_2=0,\quad\textrm{for any}\,\, t,
\end{equation}
and giving uniqueness.\\

Existence of $z^{\epsilon}$ is a straightforward of adaptations of results of Ladyzenskaja, Sollonnikov and Ural' Ceva \cite{LSU} or Lions \cite{J.L.L} on a time interval of length  $\epsilon.$\\
Then, let us consider  the function $Z^{\epsilon}=Z^{\epsilon}(t,x)=\mathcal{S}(t,\frac{t}{\epsilon},x)$ where $\mathcal{S}$ is solution to (\ref{3.104.1}). We obtain \begin{equation}\frac{\partial Z^{\epsilon}}{\partial t}=\frac{\partial S}{\partial t}(t,\frac{t}{\epsilon},x)+\frac{1}{\epsilon}\frac{\partial S}{\partial\theta}(t,\frac{t}{\epsilon},x)\end{equation} Using equation (\ref{3.104}) we deduce that $Z^{\epsilon}$ is solution to
\begin{equation}\frac{\partial Z^{\epsilon}}{\partial t}-\frac{1}{\epsilon^{2}}\nabla\cdot\Big(\tilde{\mathcal{A}}_{\epsilon}\nabla Z^{\epsilon}\Big)= \frac{1}{\epsilon^{2}}\nabla\cdot\tilde{\mathcal{C}}_{\epsilon}+\frac{\partial S}{\partial t}\end{equation}
then we deduce that 
\begin{equation}\label{zZ}\left\{\begin{array}{ccc}\ds \frac{\partial(z^{\epsilon}- Z^{\epsilon})}{\partial t}-\frac{1}{\epsilon^{2}}\nabla\cdot\Big(\tilde{\mathcal{A}}_{\epsilon}\nabla(z^{\epsilon}- Z^{\epsilon})\Big)= \frac{\partial S}{\partial t}\\
(z^{\epsilon}-Z^{\epsilon})_{|t=0}=z_{0}-\mathcal{S}(0,0,x).\end{array}\right.\end{equation}
Multiplying (\ref{zZ}) by $z^{\epsilon}- Z^{\epsilon}$ and integrating over $\torus^{2},$ we have
\begin{equation}\label{PM}\frac{d\|z^{\epsilon}- Z^{\epsilon}\|_{2}^{2}}{dt}+\frac{1}{\epsilon^{2}}\int_{\torus^{2}}\tilde{\mathcal{A}}_{\epsilon}|\nabla(z^{\epsilon}- Z^{\epsilon}|^{2}dx=\int_{\torus^{2}}\frac{\partial\mathcal{S}}{\partial t}(z^{\epsilon}- Z^{\epsilon})dx\end{equation}
which gives
\begin{equation}
\frac{d\|z^{\epsilon}- Z^{\epsilon}\|_{2}^{2}}{dt}\leq \sqrt{\epsilon\Big(\frac{\gamma+\epsilon\gamma^{3}}{\sqrt{\widetilde G_{thr}}}+(\gamma^2+\epsilon^2\gamma^4)\Big)} \|z^{\epsilon}- Z^{\epsilon}\|_{2}.\end{equation} Then we have
\begin{equation}\|z^{\epsilon}(t,\cdot)- Z^{\epsilon}(t,\cdot)\|_{2}^{2}\leq \|z_{0}-\mathcal{S}(0,0,x)\|_{2}\sqrt{\epsilon\Big(\frac{\gamma+\epsilon\gamma^{3}}{\sqrt{\widetilde G_{thr}}}+(\gamma^2+\epsilon^2\gamma^4)\Big)}T.\end{equation}
As $\|\mathcal{S}\|_{L^{\infty}_{\#}(\mathbb{R},L^{2}(\torus^{2}))}\leq\frac{\gamma}{\sqrt{\widetilde{G}_{thr}}}$ when $\epsilon\rightarrow0,$ then (\ref{3.1080}) is true.
\end{proof}
\section{Homogenization for long term dynamics of dunes, proof of theorem \,\ref{thAsyBeh1}}
We consider  equation 
(\ref{L1})
where $\mathcal{A}^{\epsilon}$ and $\mathcal{C}^{\epsilon}$ are defined by formulas (\ref{La1}) coupled with (\ref{La3}) and (\ref{La4}) coupled with (\ref{La5}).
Our aim consists in deducing the equations satisfied by the limit of $z^{\epsilon}$ solution to  (\ref{L1}) as $\epsilon\longrightarrow0.$
\\~\\
It is obvious that 
\begin{eqnarray}\label{H2}
{\mathcal{A}^{\epsilon}(t,x)\,\,\textrm{two-scale converges to}\,\,  \widetilde{\mathcal{A}}(t,\theta,x)\in L^{\infty}([0,T],L^{\infty}_\#(\R,L^2(\torus^{2}))) {} }~~~~~~~~
    \nonumber\\
\textrm{and}\,\, \mathcal{C}^{\epsilon}(t,x)\,\,\textrm{two-scale converges to}\,\, \widetilde{\mathcal{C}}(t,\theta,x),
\end{eqnarray}
with 
\begin{equation}\label{EqHom}
\widetilde{\mathcal{A}}(t,\theta,x)=a\,g_a(|\mathcal{U}(t,\theta,x)|) \,\, \textrm{and}\,\,\widetilde{\mathcal{C}}(t,\theta,x)=c\,g_c(|\mathcal{U}(t,\theta,x)|)\,\frac{\mathcal{U}(t,\theta,x)}{|\mathcal{U}(t,\theta,x)|}.
\end{equation}
Assumptions (\ref{suppl3}) and (\ref{suppl3.1}) have the following equivalence here:
\begin{equation}\label{suppl3hom}
\Theta=[0, T)\times\{\theta\in\mathbb{R}:\,\,\widetilde{\mathcal{A}}(\cdot,\theta,\cdot)=0\}\times\torus^2,
\end{equation}and
\begin{equation}\label{suppl3.1hom}
\Theta_{thr}=\{(t,\theta,x)\in [0,T)\times\mathbb{R}\times\torus^2\,\,\text{such that}\,\,\widetilde{\mathcal{A}}(t,\theta,x)<\widetilde{G}_{thr}\}.
\end{equation}

Moreover, we notice that because of (\ref{suppl1}) and (\ref{suppl2})
\begin{equation}\label{suppl3bis}
\widetilde{\mathcal{A}}(t,\theta,x)=0\,\,\text{if and only if}\,\,(t,\theta,x)\,\,\in\,\,\Theta.
\end{equation}
 We have the following theorem.
\begin{theorem}\label{thHomSecHom1.}
Under assumptions (\ref{Hy1}), (\ref{Hy2}), (\ref{Hy3}), (\ref{H2}), (\ref{EqHom}) and (\ref{suppl3bis}),  
for any $T,$ not depending on $\epsilon,$  the sequence $(z^{\epsilon})$ of solutions to (\ref{L1}), with coefficients given by (\ref{La1}) coupled with (\ref{La3}) and (\ref{La4}) coupled with (\ref{La5}), two-scale converges to the profile $U\in L^{\infty}([0,T],L^{\infty}_\#(\R,L^2(\torus^{2})))$ solution to
\begin{equation}
\label{ee179}
-\nabla\cdot(\widetilde{\mathcal{A}}\nabla U)=\nabla \cdot\widetilde{\mathcal{C}}\quad \text{on}\,\,\Big([0,T)\times \mathbb{R}\times\torus^2\Big)\backslash\Theta,\end{equation}
\begin{equation}\frac{\partial U}{\partial\theta}=0\quad \text{on}\,\,\Theta_{thr},\end{equation}
\begin{equation}\int_0^1\int_{\torus^2}Ud\theta\,dx\,\,=\,\,\int_{\torus^2}z_0dx,
\end{equation} 
where $\widetilde{\mathcal{A}}$ and $\widetilde{\mathcal{C}}$ are given by (\ref{EqHom}); $\Theta$ and $\Theta_{thr}$ are given by (\ref{suppl3hom}) and (\ref{suppl3.1hom}).
\end{theorem}
\begin{proof}.\,(of Theorem\,\ref{thAsyBeh1}). Theorem \ref{thAsyBeh1} is a direct consequence of theorem\,\ref{thHomSecHom1.}
\end{proof}
\begin{proof}.\,(of Theorem\,\ref{thHomSecHom1.}). Multiplying  (\ref{L1}) by $\psi^{\epsilon}(t,x)=\psi(t,\frac{t}{\epsilon},x)$ regular with compact support in $[0,T)\times\torus^{2}$ and 1-periodic in  $\theta,$ we obtain 
\begin{equation}\label{W.F}-\int_{\torus^{2}}z_{0}(x)\psi(0,0,x)dx-\int_{\torus^{2}}\int_{0}^{T}\frac{\partial\psi^{\epsilon}}{\partial t}z^{\epsilon}dt\,dx+$$
$$\frac{1}{\epsilon^{2}}\int_{\torus^{2}}\int_{0}^{T}\mathcal{A}_{\epsilon}\nabla z^{\epsilon}\nabla\psi^{\epsilon}dt\,dx=\frac{1}{\epsilon^{2}}\int_{\torus^{2}}\int_{0}^{T}\Big(\nabla\cdot\mathcal{C}_{\epsilon}\Big)\psi^{\epsilon}dx.
\end{equation}
Using the Green formula and \begin{equation}\label{H5}
 \frac{\partial \psi^{\epsilon}}{\partial t}
 =\left(\frac{\partial \psi}{\partial t}\right)^{\epsilon}
 +\frac{1}{\epsilon}\left(\frac{\partial\psi}{\partial \theta}\right)^{\epsilon},
\end{equation}
where
\begin{gather}\label{H4a}
   \left(\frac{\partial\psi}{\partial t}\right)^{\epsilon}(t,x)
   =\frac{\partial\psi}{\partial t}(t,\frac{t}{\epsilon},x)\,\,\textrm{ and }
     \left(\frac{\partial\psi}{\partial\theta}\right)^{\epsilon}(t,x)=\frac{\partial\psi}{\partial\theta}(t,\frac{t}{\epsilon},x),
\end{gather}
we obtain
\begin{eqnarray}\label{H6}
    \lefteqn{ \int_{\torus^{2}}\int_{0}^{T}\Big(\Big(\frac{\partial\psi}{\partial t}\Big)^{\epsilon}+\frac{1}{\epsilon}\Big(\frac{\partial\psi}{\partial \theta}\Big)\Big)z^{\epsilon}\,dt\,dx+\frac{1}{\epsilon^{2}}\int_{\torus^{2}}\int_{0}^{T}z^{\epsilon}\nabla\cdot\Big(\mathcal{A}_{\epsilon}\nabla\psi^{\epsilon}\Big)dt\,dx  {}}~~~~~~~~
    \nonumber\\
    & &+
\frac{1}{\epsilon^{2}}\int_{\torus^{2}}\int_{0}^{T}\Big(\nabla\cdot\mathcal{C}_{\epsilon}\Big)\psi^{\epsilon}dt\,dx=-\int_{\torus^{2}}z_{0}(x)\psi(0,0,x)dx\end{eqnarray}
Multiplying by $\epsilon^{2}$ and using the two-scale convergence  due to Nguetseng \cite{nguetseng:1989}, Allaire\,\cite{allaire:1992}, Fr\'enod, Raviart and Sonnendrucker \cite{FRS:1999}, as $z^{\epsilon}$ is bounded in $L^{\infty}(0,T,L^{2}(\torus^{2}))$, there exists a profile $U(t,\theta,x)$, periodic of period 1 with respect to $\theta$, such that for all 
$\psi(t,\theta,x),$ regular with compact support with respect to $(t,x)$ and  periodic of period 1 
with respect to $\theta$, we have
\begin{equation}-\int_{\torus^{2}}\int_{0}^{T}\int_{0}^{1}U\nabla\cdot\Big(\tilde{\mathcal{A}}\nabla\psi\Big)d\theta\, dt\,dx=\int_{\torus^{2}}\int_{0}^{T}\int_{0}^{1}\big(\nabla\cdot\tilde{\mathcal{C}}\big)\psi\,\,d\theta\,dt\,dx,\end{equation} then
\begin{equation}\label{ee179bis}-\nabla\cdot\Big(\tilde{\mathcal{A}}\nabla U\Big)=\nabla\cdot\tilde{\mathcal{C}},\end{equation}
with $\tilde{\mathcal{A}}$ and $\tilde{\mathcal{C}}$ given by (\ref{EqHom}).\\
Since $\widetilde{\mathcal{A}}$ and $\widetilde{\mathcal{C}}$ vanish on $\Theta,$ we deduce (\ref{ee179}) from (\ref{ee179bis}).\\
Moreover, because of (\ref{Hy2}), in points where  $\widetilde{\mathcal{A}}(t,\theta,x)<\widetilde{G}_{thr},\quad \nabla\cdot\widetilde{\mathcal{C}}=0$ and $\widetilde{\mathcal{A}}$ does not depend on $t$ and $x.$ Hence $U$ depends only on $\theta.$ In other words, \begin{equation}\label{suppl4.1}U(t,\theta,x)=U(\theta)\,\,\text{on}\,\,\Theta_{thr}.\end{equation}
Taking now test functions $\psi$ not depending on $x$ in (\ref{W.F}), the two last terms of the left hand side of (\ref{H6}) vanish. Then passing to the limit, we obtain the weak formulation of
\begin{equation}\label{suppl5}
\frac{\partial\Big(\int_{\torus^2}U(t,\theta,x)dx\Big)}{\partial\theta}=0
\end{equation} which yields because of (\ref{suppl4.1})
\begin{equation}\label{suppl6}
\frac{\partial U}{\partial \theta}=0\,\,\text{on}\,\,\Theta_{thr}.
\end{equation}
Finally, taking test function $\psi$ depending only on $t$ we obtain
\begin{equation}\label{suppl7}
\int_0^1\int_{\torus^2}U(t,\theta,x)d\theta\,dx=\int_{\torus^2}z_0(x)dx,
\end{equation}
ending the proof of the theorem.
\end{proof}
\section{Homogenization and corrector result for mean-term
 dynamics of dunes, proof of theorem\,\ref{thHomSecHom1} and \ref{th3.1}}
Making the same as in the begining of section \ref{secExEs1}, setting: \begin{equation}\label{L2}
    \mathcal{A}^{\epsilon}(t,x)=\widetilde{\mathcal{A}}_{\epsilon}(t,\frac{t}{\sqrt{\epsilon}},\frac{t}{\epsilon},x),
\end{equation}
and \begin{equation}\label{L4}
    \mathcal{C}^{\epsilon}(t,x)=\widetilde{\mathcal{C}}_{\epsilon}(t,\frac{t}{\sqrt{\epsilon}},\frac{t}{\epsilon},x),
\end{equation}
where
\begin{equation}\label{L3}\widetilde{\mathcal{A}}_{\epsilon}(t,\tau,\theta,x)=a(1-b\sqrt{\epsilon}\mathcal{M}(t,\tau,\theta,x))\,
g_a(|\mathcal{U}(t,\tau,\theta,x)|),
\end{equation} 
and
\begin{equation}\label{L5} \widetilde{\mathcal{C}}_{\epsilon}(t,\tau,\theta,x)=
c(1-b\sqrt{\epsilon}\mathcal{M}(t,\tau,\theta,x))\,g_c(|\mathcal{U}(t,\tau,\theta,x)|)\,
\frac{\mathcal{U}(t,\tau,\theta,x)}{|\mathcal{U}(t,\tau,\theta,x)|},
\end{equation}
equation (\ref{eq5}) with initial condition (\ref{eq11}) can be set in the form
\begin{equation}\label{L01}\left\{\begin{array}{cc}
    \ds \frac{\partial z^{\epsilon}}{\partial t}-\frac{1}{\epsilon}\nabla\cdot(\mathcal{A}^{\epsilon}\nabla z^{\epsilon})=\frac{1}{\epsilon}\nabla\cdot\mathcal{C}^{\epsilon},\\
    \ds z^{\epsilon}_{|t=0}=z_{0}.
\end {array}\right.\end{equation}
Under assumptions (\ref{eq2}) and (\ref{eq8}),  $\widetilde{\mathcal{A}}_{\epsilon}$ and $\widetilde{\mathcal{C}}_{\epsilon}$
given by (\ref{L3}) and (\ref{L5})  satisfy the following hypotheses:
\begin{equation}\label{Hy01}\left\{\begin{array}{ccc}
\tau\longmapsto(\widetilde{\mathcal{A}}_{\epsilon},\widetilde{\mathcal{C}}_{\epsilon})\,\,\textrm{is periodic of period}\,\,1,\\
\theta\longmapsto(\widetilde{\mathcal{A}}_{\epsilon},\widetilde{\mathcal{C}}_{\epsilon})\,\,\textrm{is periodic of period}\,\,1,
\\

x\longmapsto(\widetilde{\mathcal{A}}_{\epsilon},\widetilde{\mathcal{C}}_{\epsilon})\,\,\textrm{defined on}\,\,
\torus^2,\,\,\\

\ds|\widetilde{\mathcal{A}}_{\epsilon}|\leq\gamma,\,\,|\widetilde{\mathcal{C}}_{\epsilon}|\leq\gamma,\,\,\left|\frac{\partial\widetilde{\mathcal{A}}_{\epsilon}}{\partial t}\right|\leq\gamma,\,\,\left|\frac{\partial\widetilde{\mathcal{C}}_{\epsilon}}{\partial t}\right|\leq\gamma,\left|\frac{\partial\nabla\widetilde{\mathcal{A}}_{\epsilon}}{\partial t}\right|\leq\gamma.\\

\ds\left|\frac{\partial\widetilde{\mathcal{A}}_{\epsilon}}{\partial \theta}\right|\leq\gamma,\,\,
    \left|\frac{\partial\widetilde{\mathcal{C}}_{\epsilon}}{\partial \theta}\right|\leq\gamma,\,\,|\nabla\widetilde{\mathcal{A}}_{\epsilon}|\leq\gamma,\,\,|\nabla\cdot\widetilde{\mathcal{C}}_{\epsilon}|\leq\gamma,\,\,\left|\frac{\partial\nabla\cdot\widetilde{\mathcal{C}}_{\epsilon}}{\partial t}\right|\leq\gamma,
\end{array}\right.\end{equation}
\begin{equation}\label{Hy02}\left\{\begin{array}{ccc}
\exists \widetilde{G}_{thr}, ~ \theta_{\alpha}<\theta_{\omega}\in [0,1]\,\,\textrm{such that}\,\, 
\forall\,\,\theta\in [\theta_{\alpha},\theta_{\omega}]\Longrightarrow
 \widetilde{\mathcal{A}}_{\epsilon}(t,\tau,\theta,x)\geq\widetilde{G}_{thr},
 \\
 \ds
 \widetilde{\mathcal{A}}_{\epsilon}(t,\tau,\theta,x)\leq\widetilde{G}_{thr}\Longrightarrow\left\{\begin{array}{ccc}
\ds\frac{\partial\widetilde{\mathcal{A}}_{\epsilon}}{\partial t}(t,\tau,\theta,x)=\ds\frac{\partial\widetilde{\mathcal{A}}_{\epsilon}}{\partial \tau}(t,\tau,\theta,x)=0,\,\,
\nabla\widetilde{\mathcal{A}}_{\epsilon}(t,\tau,\theta,x)=0, \vspace{3pt}\\
\ds\frac{\partial\widetilde{\mathcal{C}}_{\epsilon}}{\partial t}(t,\tau,\theta,x)=\ds\frac{\partial\widetilde{\mathcal{C}}_{\epsilon}}{\partial \tau}(t,\tau,\theta,x)=0,\,\,
\nabla\cdot\widetilde{\mathcal{C}}_{\epsilon}(t,\tau,\theta,x)=0 ,\\
\end{array}\right.\end{array}\right.\end{equation}
 \begin{equation}\label{Hy03}\left\{\begin{array}{ccc}
\ds
 |\widetilde{\mathcal{C}}_{\epsilon}|\leq\gamma|\widetilde{\mathcal{A}}_{\epsilon}|, ~
 |\widetilde{\mathcal{C}}_{\epsilon}|^{2}\leq\gamma|\widetilde{\mathcal{A}}_{\epsilon}|,~
 |\nabla\widetilde{\mathcal{A}}_{\epsilon}|\leq\gamma|\widetilde{\mathcal{A}}_{\epsilon}|,~
 \Big|\frac{\partial\widetilde{\mathcal{A}}_{\epsilon}}{\partial t}\Big|\leq\gamma|\widetilde{\mathcal{A}}_{\epsilon}|,
\\
\ds
 \Big|\frac{\partial(\nabla\widetilde{\mathcal{A}}_{\epsilon})}{\partial t}
       \Big|^{2}\leq\gamma|\widetilde{\mathcal{A}}_{\epsilon}|,~
 \Big|\nabla \cdot\widetilde{\mathcal{C}}_{\epsilon}\Big|\leq\gamma|\widetilde{\mathcal{A}}_{\epsilon}|,~
 \Big|\frac{\partial\widetilde{\mathcal{C}}_{\epsilon}}{\partial t}\Big|\leq\gamma|\widetilde{\mathcal{A}}_{\epsilon}|,~
 \Big|\frac{\partial\widetilde{\mathcal{C}}_{\epsilon}}{\partial t}\Big|^{2}\leq\gamma^2|\widetilde{\mathcal{A}}_{\epsilon}|\\
\Big|\ds \frac{\partial\widetilde{\mathcal{A}}_{\epsilon}}{\partial\tau}\Big|^{2}\leq\epsilon\gamma|\widetilde{\mathcal{A}}_{\epsilon}|,\,\,\ds \Big|\frac{\partial\nabla\widetilde{\mathcal{A}}_{\epsilon}}{\partial\tau}\Big|^{2}\leq\epsilon\gamma|\widetilde{\mathcal{A}}_{\epsilon}|.
\end{array}\right.\end{equation}
For (\ref{L01}), if hypotheses (\ref{Hy01}), (\ref{Hy02}) and (\ref{Hy03}) are satisfied, an existence and uniqueness result is given in \cite{FaFreSe}.
\subsection{Homogenization }
Let us consider equation (\ref{L01}) with  $\mathcal{A}_{\epsilon}$ and $\mathcal{C}_{\epsilon}$ given by (\ref{L2}) and (\ref{L4});
\begin{equation}\label{H20a}
\theta\longmapsto\,\,\widetilde{\mathcal{A}},\,\,\widetilde{\mathcal{C}}\,\,\text{is periodic of period 1},$$
$$\tau\longmapsto\,\,\widetilde{\mathcal{A}},\,\,\widetilde{\mathcal{C}}\,\,\text{is periodic of period 1}, 
\end{equation}
\begin{eqnarray}\label{H2a}
{\mathcal{A}^{\epsilon}(t,x)\,\,\textrm{3-scale converges to}\,\,  \widetilde{\mathcal{A}}(t,\tau,\theta,x)\in L^{\infty}([0,T]\times\mathbb{R},L^{\infty}_\#(\R,L^2(\torus^{2}))) {} }~~~~~~~~
    \nonumber\\
\textrm{and}\,\, \mathcal{C}^{\epsilon}(t,x)\,\,\textrm{3-scale converges to}\,\, \widetilde{\mathcal{C}}(t,\tau,\theta,x),
\end{eqnarray}
with 
\begin{equation}\label{EqHoma}
\widetilde{\mathcal{A}}(t,\tau,\theta,x)=a\,g_a(|\mathcal{U}(t,\tau,\theta,x)|) \,\, \textrm{and}\,\,\widetilde{\mathcal{C}}(t,\tau,\theta,x)=c\,g_c(|\mathcal{U}(t,\tau,\theta,x)|)\,\frac{\mathcal{U}(t,\tau,\theta,x)}{|\mathcal{U}(t,\tau,\theta,x)|}.
\end{equation}

\begin{theorem}\label{thHomSecHom12}
Under assumptions (\ref{Hy01}), (\ref{Hy02}), (\ref{Hy03}),(\ref{H20a}), (\ref{H2a}) and (\ref{EqHoma}), 
for any $T,$ not depending on $\epsilon,$  the sequence $(z^{\epsilon})$ of solutions to (\ref{L01}), with coefficients given by (\ref{L2}) coupled with (\ref{L3}) and (\ref{L4}) coupled with (\ref{L5}), 3-scale converges to the profile $U\in L^{\infty}([0,T]\times\mathbb{R},L^{\infty}_\#(\R,L^2(\torus^{2})))$ solution to
\begin{equation}
\label{ee179tri}
\frac{\partial U}{\partial\theta}-\nabla\cdot(\widetilde{\mathcal{A}}\nabla U)=\nabla \cdot\widetilde{\mathcal{C}},
\end{equation} where $\widetilde{\mathcal{A}}$ and $\widetilde{\mathcal{C}}$ are given by (\ref{EqHoma}).
\end{theorem}
\begin{proof}.\,(of Theorem\,\ref{thHomSecHom1}). Theorem \ref{thHomSecHom1} is a direct consequence of theorem\,\ref{thHomSecHom12}.
\end{proof}
\begin{proof}.\,(of Theorem\,\ref{thHomSecHom12}).
Considering test functions $\psi^{\epsilon}(t,x)=\psi(t,\frac{t}{\sqrt{\epsilon}},\frac{t}{\epsilon},x)$ for all $\psi(t,\tau,\theta,x)$ regular with compact support on  $[0,T)\times\torus^{2}$ and periodic of period $1$ with respect to  $\tau$ and $\theta.$
\begin{equation}\label{hye}\frac{\partial\psi^{\epsilon}}{\partial t}=\big(\frac{\partial\psi}{\partial t}\big)^{\epsilon}+\frac{1}{\sqrt{\epsilon}}\big(\frac{\partial\psi}{\partial\tau}\big)^{\epsilon}+\frac{1}{\epsilon}\big(\frac{\partial\psi}{\partial\theta}\big)^{\epsilon},
\end{equation}
where \begin{equation}
\big(\frac{\partial\psi}{\partial t}\big)^{\epsilon}(t,x)=\frac{\partial\psi}{\partial t}(t,\frac{t}{\sqrt{\epsilon}},\frac{t}{\epsilon},x),\qquad\big(\frac{\partial\psi}{\partial\tau}\big)^{\epsilon}=\frac{\partial\psi}{\partial\tau}(t,\frac{t}{\sqrt{\epsilon}},\frac{t}{\epsilon},x),\qquad
\big(\frac{\partial\psi}{\partial\theta}\big)^{\epsilon}=\frac{\partial\psi}{\partial\theta}(t,\frac{t}{\sqrt{\epsilon}},\frac{t}{\epsilon},x).
\end{equation}
Multiplying (\ref{L01}) by $\psi^{\epsilon}(t,\frac{t}{\sqrt{\epsilon}},\frac{t}{\epsilon},x)$ and integrating on $[0,T)\times\torus^{2},$ we get
\begin{gather*}-\int_{\torus^{2}}z_{0}(x)\psi(0,0,0,x)dx-\int_{\torus^{2}}\int_{0}^{T}\frac{\partial\psi^{\epsilon}}{\partial t}z^{\epsilon}dt\,dx-
\frac{1}{\epsilon}\int_{\torus^{2}}\int_{0}^{T} z^{\epsilon}\nabla\cdot\Big(\mathcal{A}^{\epsilon}\nabla\psi^{\epsilon}\Big)dt\,dx\\
=\frac{1}{\epsilon}\int_{\torus^{2}}\int_{0}^{T}\nabla\cdot\mathcal{C}_{\epsilon}\psi^{\epsilon}dt\,dx.
\end{gather*}
Replacing  $\frac{\partial\psi^{\epsilon}}{\partial t}$ by the relation (\ref{hye}), we have
\begin{gather*}\int_{\torus^{2}}\int_{0}^{T}z^{\epsilon}\Big[\big(\frac{\partial\psi}{\partial t}\big)^{\epsilon}+\frac{1}{\sqrt{\epsilon}}\big(\frac{\partial\psi}{\partial\tau}\big)^{\epsilon}+\frac{1}{\epsilon}\big(\frac{\partial\psi}{\partial\theta}\big)^{\epsilon}\Big]dt\,dx+\frac{1}{\epsilon}\int_{\torus^{2}}\int_{0}^{T}z^{\epsilon}\nabla\cdot\Big(\mathcal{A}^{\epsilon}\nabla\psi^{\epsilon}\Big)dt\,dx\\
+\frac{1}{\epsilon}\int_{\torus^{2}}\int_{0}^{T}\nabla\cdot\mathcal{C}^{\epsilon}\psi^{\epsilon}dt\,dx=-\int_{\torus^{2}}z_{0}(x)\psi(0,0,0,x)dx.
\end{gather*}
Multiplying by  $\epsilon$ we have
\begin{gather*}\int_{\torus^{2}}\int_{0}^{T}z^{\epsilon}\Big[\epsilon\big(\frac{\partial\psi}{\partial t}\big)^{\epsilon}+\sqrt{\epsilon}\big(\frac{\partial\psi}{\partial\tau}\big)^{\epsilon}+\big(\frac{\partial\psi}{\partial\theta}\big)^{\epsilon}+ \nabla\cdot\Big(\mathcal{A}^{\epsilon}\nabla\psi^{\epsilon}\Big)\Big]dt\,dx+\\
\int_{\torus^{2}}\int_{0}^{T}\nabla\cdot\mathcal{C}^{\epsilon}\psi^{\epsilon}dt\,dx=-\epsilon\int_{\torus^{2}}z_{0}(x)\psi(0,0,0,x)dx.
\end{gather*}
The functions $\big(\frac{\partial\psi}{\partial t}\big)^{\epsilon},\,\,\big(\frac{\partial\psi}{\partial\tau}\big)^{\epsilon}\,\,\textrm{and}\,\,\big(\frac{\partial\psi}{\partial\theta}\big)^{\epsilon}$ are periodic with respect to the two variables $\tau,\,\,\theta.$ Here we use the  $3$-scales convergence see \cite{J.L.D.L.L.P}.\\
Taking the limit as $\epsilon\rightarrow0,$ using the $3$-scales convergence, we have 
$$\int_{\torus^{2}}\int_{0}^{T}\int_{[0,1]^{2}}\Big(U\frac{\partial\psi}{\partial\theta}+U\nabla\cdot\big(\tilde{\mathcal{A}}\nabla\psi\big)\Big)d\tau\, d\theta\, dt\,dx=\int_{\torus^{2}}\int_{0}^{T}\int_{[0,1]^{2}}\tilde{\mathcal{C}}\cdot\nabla\psi\, d\tau\, d\theta\, dt\,dx.$$ Then, the limit $U$ of  $z^{\epsilon}$ solution to (\ref{L1}) satisfies the following equation
\begin{equation}\label{ff}
\frac{\partial U}{\partial \theta}-\nabla\cdot\Big(\tilde{\mathcal{A}}\nabla U\Big)=\nabla\cdot\tilde{\mathcal{C}}.\end{equation}
\end{proof}
There is indeed existence and uniqueness of the equation (\ref{ff}) according to the application of the theorem 3.15 of \cite{FaFreSe}; thus (\ref{ff}) is  the homogenized equation. In (\ref{ff}), $\tau\,\,\text{and}\,\,t$ are only parameters. 
\subsection{A corrector result}
Considering equation (\ref{L01}) with coefficients (\ref{L2}) and (\ref{L4}) and hypothesis (\ref{H2a}) leads to
\begin{equation}\mathcal{A}^{\epsilon}(t,x)=\widetilde{\mathcal{A}}^{\epsilon}(t,x)+\sqrt{\epsilon}\widetilde{\mathcal{A}}^{\epsilon}_{1}(t,x)+\epsilon\widetilde{\mathcal{A}}^{\epsilon}_{2}(t,x),\end{equation}
\begin{equation}\mathcal{C}^{\epsilon}(t,x)=\widetilde{\mathcal{C}}^{\epsilon}(t,x)+\sqrt{\epsilon}\widetilde{\mathcal{C}}^{\epsilon}_{1}(t,x)+\epsilon\widetilde{\mathcal{C}}^{\epsilon}_{2}(t,x)\end{equation}
where 
\begin{equation}\label{cr3}\widetilde{\mathcal{A}}^{\epsilon}(t,x)=\widetilde{\mathcal{A}}(t,\frac{t}{\sqrt{\epsilon}},\frac{t}{\epsilon},x),\qquad\widetilde{\mathcal{C}}^{\epsilon}(t,x)=\widetilde{\mathcal{C}}(t,\frac{t}{\sqrt{\epsilon}},\frac{t}{\epsilon},x)\end{equation}
\begin{equation}\label{cr4}\widetilde{\mathcal{A}}^{\epsilon}_{1}(t,x)=\widetilde{\mathcal{A}}_{1}(t,\frac{t}{\sqrt{\epsilon}},\frac{t}{\epsilon},x),\qquad\widetilde{\mathcal{C}}^{\epsilon}_{1}(t,x)=\widetilde{\mathcal{C}}_{1}(t,\frac{t}{\sqrt{\epsilon}},\frac{t}{\epsilon},x)\end{equation}
\begin{equation}\label{cr5}\widetilde{\mathcal{A}}^{\epsilon}_{2}(t,x)=\widetilde{\mathcal{A}}_{2}(t,\frac{t}{\sqrt{\epsilon}},\frac{t}{\epsilon},x),\qquad\widetilde{\mathcal{C}}^{\epsilon}_{2}(t,x)=\widetilde{\mathcal{C}}_{2}(t,\frac{t}{\sqrt{\epsilon}},\frac{t}{\epsilon},x)\end{equation}
Because of hypotheses (\ref{Hy01}), (\ref{Hy02}) and (\ref{Hy03}),\,\,
$\widetilde{\mathcal{A}},\,\,\widetilde{\mathcal{A}}_{1},\,\,\widetilde{\mathcal{A}}_{2},\,\,\widetilde{\mathcal{A}}^{\epsilon},\,\,\widetilde{\mathcal{A}}^{\epsilon}_{1},\,\,\widetilde{\mathcal{A}}^{\epsilon}_{2},\,\,\widetilde{\mathcal{C}},\,\,\widetilde{\mathcal{C}}_{1},\,\,\widetilde{\mathcal{C}}_{2},\,\,\widetilde{\mathcal{C}}^{\epsilon},\,\,\widetilde{\mathcal{C}}^{\epsilon}_{1} \,\,\textrm{and}\,\,\widetilde{\mathcal{C}}^{\epsilon}_{2}$ are regular and bounded coefficients.
\begin{theorem}\label{th3.1bis}
Under assumptions (\ref{Hy01}), (\ref{Hy02}), (\ref{Hy03}),(\ref{H20a}), (\ref{H2a}) and (\ref{EqHoma}),
considering function $z^{\epsilon}\in
L^{\infty}([0,T),L^{2}(\torus^{2})),$ solution to (\ref{L1})
and function $U^{\epsilon}\in
L^{\infty}([0,T],L^{\infty}_\#(\R,L^2(\torus^{2})))$ defined by
$U^{\epsilon}(t,x)=U(t,\frac{t}{\sqrt{\epsilon}},\frac{t}{\epsilon},x),$
where $U$ is the solution to (\ref{ee179tri}), the following
estimate is satisfied:
\begin{equation}
\Big\|\frac{z^{\epsilon}-U^{\epsilon}}{\sqrt{\epsilon}}\Big\|_{
L^{\infty}([0,T),L^{2}(\torus^{2}))}\,\,\leq\alpha,
\end{equation}
where $\alpha$ is a constant not depending on $\epsilon.$\\
Furthermore
\begin{equation}
\frac{z^{\epsilon}-U^{\epsilon}}{\sqrt{\epsilon}}\quad\textrm{
3-scale converges to a profile}\,\,U_{\frac{1}{2}}\in
L^{\infty}([0,T]\times\mathbb{R},L^{\infty}_\#(\R,L^2(\torus^{2}))),
\end{equation}
which is the unique solution to
\begin{equation}\frac{\partial U_{\frac{1}{2}}}{\partial\theta}-\nabla\cdot\Big(\widetilde{\mathcal{A}}\nabla U_{\frac{1}{2}}\Big)=\nabla\cdot\widetilde{\mathcal{C}}^{1}+\nabla\cdot\Big(\widetilde{\mathcal{A}}_{1}\nabla U\Big)-\frac{\partial U}{\partial\tau}.
\end{equation}
\end{theorem}
\begin{proof}.\,(of Theorem\ref{th3.1}). Theorem\,\ref{th3.1} is a direct consequence of theorem\,\ref{th3.1bis}.
\end{proof}
\begin{proof}.\,(of Theorem\,\ref{th3.1bis}).
Using the relations (\ref{cr3}), (\ref{cr4}) and (\ref{cr5}), equation (\ref{L01}) becomes
\begin{equation}\label{cr6}
\frac{\partial z^{\epsilon}}{\partial t}-\frac{1}{\epsilon}\nabla\cdot\Big(\widetilde{\mathcal{A}}^{\epsilon}\nabla z^{\epsilon}\Big)=\frac{1}{\epsilon}\Big(\nabla\cdot\widetilde{\mathcal{C}}^{\epsilon}+\sqrt{\epsilon}\nabla\cdot\widetilde{\mathcal{C}}_{1}+\epsilon\nabla\cdot\widetilde{\mathcal{C}}_{2}+
\sqrt{\epsilon}\nabla\cdot\Big(\widetilde{\mathcal{A}}^{\epsilon}_{1}\nabla z^{\epsilon}\Big) +\epsilon\nabla\cdot\Big(\widetilde{\mathcal{A}}^{\epsilon}_{2}\nabla z^{\epsilon}\Big)\Big).
\end{equation}
As $U$ is solution to (\ref{ee179tri}) and taking into account that 
\begin{equation}
\frac{\partial U^{\epsilon}}{\partial t}=\Big(\frac{\partial U}{\partial t}\Big)^{\epsilon}+\frac{1}{\sqrt{\epsilon}}\Big(\frac{\partial U}{\partial\tau}\Big)^{\epsilon}+\frac{1}{\epsilon}\Big(\frac{\partial U}{\partial\theta}\Big)^{\epsilon},
\end{equation}
we obtain the following equation
\begin{equation}\label{cr8}
\frac{\partial U^{\epsilon}}{\partial t}-\frac{1}{\epsilon}\nabla\cdot\Big(\widetilde{\mathcal{A}}^{\epsilon}\nabla U^{\epsilon}\Big)=\frac{1}{\epsilon}\Big(\nabla\cdot\widetilde{\mathcal{C}}^{\epsilon}+\sqrt{\epsilon}\Big(\frac{\partial U}{\partial\tau}\Big)^{\epsilon}+\epsilon\Big(\frac{\partial U}{\partial t}\Big)^{\epsilon}\Big).
\end{equation}
Considering equation (\ref{cr6}) and (\ref{cr8}), $z^{\epsilon}-U^{\epsilon}$ is solution to
\begin{equation}\label{cr9}\frac{ \partial(\frac{z^{\epsilon}- U^{\epsilon}}{\sqrt{\epsilon}})}{\partial t}
-\frac{1}{\epsilon}\nabla\cdot\Big(\big(\widetilde{\mathcal{A}}^{\epsilon}+\sqrt{\epsilon}\widetilde{\mathcal{A}}^{\epsilon}_{1}+\epsilon\widetilde{\mathcal{A}}^{\epsilon}_{2}\big) \nabla(\frac{z^{\epsilon}- U^{\epsilon}}{\sqrt{\epsilon}})\Big)=\frac{1}{\epsilon}\Big(\nabla\cdot\widetilde{\mathcal{C}}_{1}^{\epsilon}+\sqrt{\epsilon}\nabla\cdot\widetilde{\mathcal{C}}_{2}^{\epsilon}+
\nabla\cdot\Big(\widetilde{\mathcal{A}}^{\epsilon}_{1}\nabla U^{\epsilon}\Big) +\sqrt{\epsilon}\nabla\cdot\Big(\widetilde{\mathcal{A}}^{\epsilon}_{2}\nabla U^{\epsilon}\Big)$$
$$
-\sqrt{\epsilon}\Big(\frac{\partial U}{\partial t}\Big)^{\epsilon}-\Big(\frac{\partial U}{\partial\tau}\Big)\Big).
\end{equation}
Using the fact that $U$ solution to (\ref{ee179tri}) belongs to  $L^{\infty}([0,T]\times\mathbb{R},L^{\infty}_\#(\R,L^2(\torus^{2}))),$  $U^{\epsilon}$ is solution to (\ref{cr8})  and a results of Ladyzenskaja, Solonnikov and Ural'Ceva \cite{LSU}, all the terms $\frac{\partial U}{\partial\tau},\,\,\frac{\partial U}{\partial t}\,\,$ are bounded.
The terms  $\widetilde{\mathcal{A}}^{\epsilon}_{1},\,\,\widetilde{\mathcal{A}}^{\epsilon}_{2},\,\,\widetilde{\mathcal{C}}^{\epsilon}_{1}\,\,\text{and}\,\,\widetilde{\mathcal{C}}^{\epsilon}_{2}$ are also bounded by hypotheses and then so are $\nabla\cdot\widetilde{\mathcal{C}}_{1}^{\epsilon},\,\,\nabla\cdot\widetilde{\mathcal{C}}_{2}^{\epsilon}$ and $\nabla\cdot\Big(\widetilde{\mathcal{A}}_{1}\nabla U^{\epsilon}\Big),\,\,\nabla\cdot\Big(\widetilde{\mathcal{A}}_{2}\nabla U^{\epsilon}\Big).$
Using the same arguments as in the proof of Theorem 1.1 in \cite{FaFreSe} we obtain that $\frac{z^{\epsilon}- U^{\epsilon}}{\sqrt{\epsilon}}$converges to a profile $U_{\frac{1}{2}}\in L^{\infty}([0,T]\times\mathbb{R},L^{\infty}_\#(\R,L^2(\torus^{2})))$ solution to
\begin{equation}
\frac{\partial U_{\frac{1}{2}}}{\partial\theta}-\nabla\cdot\Big(\widetilde{\mathcal{A}}\nabla U_{\frac{1}{2}}\Big)=\nabla\cdot\widetilde{\mathcal{C}}^{1}+\nabla\cdot\Big(\widetilde{\mathcal{A}}_{1}\nabla U\Big)-\frac{\partial U}{\partial\tau}.
\end{equation}
\end{proof}


\begin{thebibliography}{99}

\bibitem{allaire:1992}
G.~Allaire, \emph{{H}omogenization and two-scale convergence}, SIAM J.
  Math. Anal. \textbf{23} (1992), 1482--1518.
\bibitem{Bagnold}
R.A. Bagnold, \emph{The movement of desert sand}, Proceedings of the Royal
  Society of London A \textbf{157} (1936), 594--620.

\bibitem{FaFreSe} I.~ Faye, E.~ Fr\'enod, D.~ Seck, \emph{Singularly perturbed degenerated parabolic equations and application to seabed 
 morphodynamics in tided environment},Discrete and Continuous Dynamical Systems, Vol 29 $N^{o} 3$ March 2011, pp 1001-1030.

\bibitem{FRS:1999}
E.~Fr\'enod, P.~A.~ Raviart, and E.~Sonnendr\"ucker, \emph{Asymptotic expansion
  of the {V}lasov equation in a large external magnetic field}, J. Math. Pures
  et Appl. \textbf{80} (2001), 815--843.

\bibitem{GaddLavSw}
P.E. Gadd, W.~Lavelle, and D.J.P. Swift, \emph{Estimates of sand transport on
  the {N}ew {Y}ork shelf using near-bottom current meter observations}, J. Sed.
  Petrol. \textbf{48} (1978.), 239--252.

\bibitem{Idier}
D.~Idier, "Dunes et bancs de sables du plateau continental: observations
  in-situ et mod\'elisation num\'erique", Ph.D. thesis, INP Toulouse, France  2002.

\bibitem{IdierAsH}
D.~Idier, D.~Astruc, and S.J.M.H.~Hulcher, \emph{Influence of bed roughness on
  dune and megaripple generation}, Geophysical Research Letters \textbf{31}
  (2004), 1--5.

\bibitem{LSU}
O.~A. Ladyzenskaja, V. A.~Solonnikov, and N. N.~Ural'ceva, "Linear and
  quasi-linear equations of parabolic type", AMS Translation of
  Mathematical Monographs \textbf{23} (1968).

\bibitem{J.L.L}
J.~L.~ Lions, \emph{Remarques sur les \'equations diff\'erentielles ordinaires},
  Osaka Math. J. \textbf{15} (1963), 131--142.

\bibitem{J.L.D.L.L.P} J.~L.~ Lions, D.~Lukkassen, L.~E.~ Persson, \emph{Reiterated homogenization of monotone operators} 
C. R. Acad. Sci. Paris, t. 330, Série I, p. 675–680, 2000
Equations aux d\'eriv\'ees partielles/Partial Differential Equations


\bibitem{MePetMull}
E.~Meyer-Peter and R.~M\"uller, \emph{Formulas for bed-load transport.}, The
  Second Meeting of the International Association for Hydraulic Structures,
  Appendix 2, (1948), 39--44.

\bibitem{nguetseng:1989}
G.~Nguetseng, \emph{A general convergence result for a functional related to
  the theory of homogenization}, SIAM J. Math. Anal. \textbf{20} (1989),
  608--623.


\bibitem{Rijn1989}
L.~C. Van~Rijn, "Handbook on sediment transport by current and waves",
  Tech. Report H461:12.1--12.27, Delft Hydraulics, 1989.

\end{thebibliography}
\end{document}